\documentclass[12pt]{amsart}

\usepackage{a4wide}

\usepackage{color}

\usepackage{amssymb}
\usepackage{epsfig}
\usepackage{subfigure}

\newtheorem{thm}{Theorem}[section]
\newtheorem{dfn}[thm]{Definition}

\newtheorem{prop}[thm]{Proposition}

\theoremstyle{remark} 
\newtheorem{rmk}[thm]{Remark}

\def\cqfd{\mbox{}\nolinebreak\hfill$\Box$\medbreak\par}
\newenvironment{pf}{\noindent\textbf{Proof:}}{\cqfd}

\def\ch#1{{{\bf H}^{#1}_{\C}}}
\newcommand{\Z}{\mathbb{Z}}
\newcommand{\N}{\mathbb{N}}
\newcommand{\C}{\mathbb{C}}
\renewcommand{\H}{\mathbb{H}}
\newcommand{\R}{\mathbb{R}}
\newcommand{\Q}{\mathbb{Q}}
\newcommand{\B}{\mathcal{B}}

\title{On spherical CR uniformization of 3-manifolds}

\author{Martin Deraux}

\date{Dec 7, 2014}

\begin{document}

\begin{abstract}
  We consider the discrete representations of 3-manifold groups into
  $PU(2,1)$ that appear in the Falbel-Koseleff-Rouillier, such that the
  peripheral subgroups have cyclic unipotent holonomy. We show that
  two of these representations have conjugate images, even though they
  represent different 3-manifold groups. This illustrates the fact
  that a discrete representation $\pi_1(M)\rightarrow PU(2,1)$ with
  cyclic unipotent boundary holonomy is not in general the holonomy of
  a spherical CR uniformization of $M$.
\end{abstract}

\maketitle

\section{Introduction}

It is a difficult problem to characterize 3-manifolds which admit a
spherical CR uniformization, i.e. manifolds which can occur as the
manifold at infinity of some infinite volume complex hyperbolic
surface (or perhaps complex hyperbolic orbifold with isolated
singularities).

Of course $S^3$ as well as lens spaces trivially admit a spherical CR
uniformizations (the orbifold fundamental group of the corresponding
complex hyperbolic quotient is then a finite cyclic group). Quotients
of the 3-dimensional Heisenberg group by lattices also provide another
somewhat trivial class, since once can think of the Heisenberg group
as a model of $\partial H^2_\C\setminus\{pt\}$; this class consists of
circle bundles over a 2-torus.
More generally, it is well known that many Seifert
fibrations occur, by taking deformations of Fuchsian subgroups of
$PU(2,1)$. For small deformations, the corresponding quotients are
disk bundles over surfaces (or more generally 2-dimensional
orbifolds), yielding circle bundles as manifolds at infinity. Note
that none of the above examples are hyperbolic manifolds.

In the last decade, R. Schwartz discovered that many hyperbolic
manifolds can occur as well, see~\cite{richBook} for a nice overview
of his constructions. The starting point was the construction of a
spherical CR uniformization of the Whitehead link
complement~\cite{richWLC}, and of a \emph{closed} hyperbolic
3-manifold~\cite{schwartz4447}. Schwartz went on to produce infinitely
many examples through a somewhat delicate Dehn surgery theorem. Recent
work of Parker and Will~\cite{parkerwill} shows that the Whitehead
link complement admits at least two distinct spherical CR
uniformizations, i.e. it occurs as the manifold at infinity of two
non-isometric complex hyperbolic orbifolds.

The figure eight knot complement was shown to admit a spherical CR
uniformization~\cite{derauxfalbel}. In fact it turns out uncountably
many pairwise non-conjugate discrete subgroups of $PU(2,1)$ yield the
figure eight knot as their manifold at infinity,
see~\cite{derauxdeformfig8}.

The general question of the classification of hyperbolic 3-manifolds
which admit a spherical CR uniformization is still widely open.

A general approach to investigate the above general question was laid
down by Falbel a few years ago. He devised a computational way to
determine all the conjugacy classes of representations of the
fundamental group of any given open open 3-manifolds (with torus ends)
with unipotent boundary holonomy. The rough idea is to use an ideal
triangulation of the manifold, so that representations are
parametrized by cross-ratios of quadruples of points in
$\partial_\infty H^2_\C$, and to write compatibility equations for
these cross-ratios to yield a representation (with appropriate
boundary holonomy conditions). For more information on this,
see~\cite{falbelfigure8},~\cite{bfg} and for another parametrization
of the representation variety, see~\cite{ggz}.

For manifolds that can be built by gluing up to three tetrahedra, the
solutions to the compatibility equations can be computed using current
computer technology (and somewhat sophisticated computational
techniques for solving polynomial systems), see~\cite{fkr}. We will
refer to this list of representations as the FKR census.

It turns out there are in fact very few discrete representations in
the FKR census. Since the representations in the census are not
representative of the whole character variety (for computational
convenience, the authors list only representations where the
peripheral holonomy consists of unipotent matrices), we will not take
the trouble of giving detailed arguments that prove non-discreteness
results.

We focus on cases where we know the representation is discrete, namely
the pairs $(M,\rho)$ in Table~\ref{tab:discrete}.
\begin{table}[htbp]
\centering
\begin{tabular}{l}
    (1) \verb|m003|, $\rho=\rho_{2\_1}$\\
    (2) \verb|m004|, $\rho=\rho_{1\_1}$\\
    (3) \verb|m004|, $\rho=\rho_{3\_1}$\\
    (4) \verb|m004|, $\rho=\rho_{4\_1}$\\
    (5) \verb|m009|, $\rho=\rho_{4\_3}$\\
    (6) \verb|m015|, $\rho=\rho_{3\_3}$
\end{tabular}
\caption{The list of discrete representations in the FKR census, for
  non-compact manifolds built out of at most three
  tetrahedra.}\label{tab:discrete}
\end{table}
We suspect that these are in fact the only discrete representations
into $PU(2,1)$ in the FKR census (apart from the representations with
finite image, which do not appear in the list in~\cite{fkr}).

Note that it is not clear whether discreteness of a representation
$\rho:\pi_1(M)\rightarrow PU(2,1)$ (even together with boundary
parabolicity) suffices to guarantee that $\rho$ is the holonomy of a
spherical CR structure on $M$. The main problem is that there is no
natural way to extend quadruples of points to
full-fledged tetrahedra in $\partial H^2_\C$. For
instance, the attempts in~\cite{falbelfigure8} and~\cite{falbelwang}
yield branched CR structures, and it is not known whether these
representations are the holonomy of an unbranched CR structure.

In fact we will be more restrictive and require not only that $\rho$
be the holonomy of a structure, but that it produce a spherical CR
{\bf uniformization} of $M$ (Schwartz call these \emph{complete}
spherical CR structures, see~\cite{richBook} for instance). We briefly
recall basic definitions pertaining to this notion.

Recall that a discrete group $\Gamma\subset PU(2,1)$, acts properly on
the complex hyperbolic plane $\ch 2$. The action extends to its
boundary at infinity $\partial_\infty \ch 2$, but it is usually no
longer proper.
\begin{dfn}
  The domain of discontinuity $\Omega_\Gamma$ is the largest open
  subset of $\partial_\infty \ch 2$ where the action is proper. Its
  complement $\Lambda_\Gamma=\partial_\infty \ch 2 - \Omega_\Gamma$ is
  called the limit set of $\Gamma$.
\end{dfn}

When the action of $\Gamma$ on $\Omega_\Gamma$ has no fixed points,
the quotient $\Gamma\setminus \Omega_\Gamma$ is a manifold, which of
course carries a CR structure inherited from the standard CR structure
on $\partial_\infty \ch 2\simeq S^3$.

\begin{dfn}
  Let $\Gamma\subset PU(2,1)$ be a discrete group whose action on $\ch
  2$ has only isolated fixed points. Then the quotient
  $\Gamma\setminus \Omega_\Gamma$ is called the {\bf manifold at
    infinity} of $\Gamma$.
\end{dfn}
We will sometimes call $\Gamma\setminus \Omega_\Gamma$ the manifold at
infinity of $\Gamma$, rather than the manifold at infinity of
$\Gamma\setminus H^2_\C$.

\begin{dfn}
  Let $\rho:\pi_1(M)\rightarrow PU(2,1)$ be a representation. We say
  that $\rho$ gives a spherical CR uniformization of $M$ if
  $\Gamma=\textrm{Im}(\rho)$ is discrete, all its fixed points in $\ch
  2$ are isolated, and the manifold at infinity of $\Gamma$ is
  homeomorphic to $M$.
\end{dfn}

We now summarize what is known about the representations that appear
in Table~\ref{tab:discrete}, in historical order
(the numbers opening each paragraph are given to follow the notation
in the table).

(2), (3), (4) The representations of $\pi_1(\verb|m004|)$ were studied
in~\cite{falbelfigure8} and~\cite{derauxfalbel}. They are all
discrete, non-faithful representations. 

(3) The group $Im(\rho_{3\_1})$ can be checked to be a
(non-elementary) normal subgroup of the Eisenstein-Picard lattice,
i.e. a normal subgroup of $PU(2,1,\Z[\omega])$
(see~\cite{falbelfigure8}). It follows that its limit set is all of
$S^3$, or equivalently that it has empty domain of discontinuity. This
makes it obvious that $\rho_{2\_1}$ is not the holonomy of a
uniformization (but it is not known whether it is the holonomy of a
spherical CR structure).

(1) A similar property holds for the image $Im(\rho_{2\_1})$ of the only
representation of $\pi_1(\verb|m003|)$ in the FKR census, which is a
normal subgroup of a lattice sometimes referred to as the sister of
the Eisenstein-Picard lattice (the Eisenstein-Picard lattice and its
sister have the same covolume, and up to conjugation they are the only
non-uniform arithmetic lattices with that volume, see~\cite{stover}).

(2), (4) The two other representations of $\pi_1(\verb|m004|)$ were
studied in~\cite{derauxfalbel}, where the author and Falbel gave a
proof that they both give a spherical CR uniformization of the figure
eight knot complement (in fact these two representations differ by
precomposition by an orientation reversing outer automorphism of
$\pi_1(\verb|m004|)$).

A more enlightening fundamental domain for the action of the group
$Im(\rho_{1\_1})$ can be obtained quite naturally as a Ford domain
centered at the fixed point of the image of a peripheral subgroup, as
worked out in~\cite{derauxdeformfig8}. That domain exhibits an
explicit horotube structure for the group, as defined
in~\cite{richBook}. Moreover, the combinatorial structure of the Ford
fundamental domain exhibits striking similarities with the structure of
the Ford domain in $H^3_\R$ for the holonomy group of the unique
complete hyperbolic structure on the figure eight knot complement.

(5) Performing the same analysis for other discrete groups with cyclic
unipotent holonomy in the FKR census, one gets the following:
\begin{thm} \label{th1}
  The representation $\rho_{4\_3}$ is the holonomy of a spherical CR
  uniformization of $M=\verb|m009|$. 
\end{thm}
Rather than using the group of the FKR census, we will give a triangle
group interpretation of that group, in the spirit of Schwartz's
constructions. This will make it easier for the reader to get
explicit matrices for the group.

The author went through the same analysis as
in~\cite{derauxdeformfig8} and worked out the combinatorics of the
Ford domain for $Im(\rho_{4\_3})$ (really, he instructed the computer
to work this out for him). The group presentation then
comes for free from the Ford domain. From that Ford domain, it is
relatively easy to compute the fundamental group of the manifold at
infinity, and to find an explicit isomorphism with the fundamental
group of \verb|m009| (of course one also needs to check that the
peripheral subgroups correspond in this isomorphism).

The details are quite unpleasant to write in a paper
(see~\cite{derauxdeformfig8} for similar arguments), so
we will not give the details of the proof of Theorem~\ref{th1}. We
only focus on proving discreteness and getting a group presentation
for $Im(\rho_{4\_3})$ (see Theorem~\ref{thm:poincare}).

Note that, unlike the case of the figure eight knot complement, the
boundary of the complex hyperbolic Ford domain does not have exactly
the same local combinatorial structure as the Ford domain of the real
hyperbolic structure (compare the 2-faces of Figure~\ref{fig:prisms}
with the shaded 2-faces in Figures~\ref{fig:faces-1},~\ref{fig:faces-2}).

(5)=(6) In view of the main results
of~\cite{schwartz4447},~\cite{derauxfalbel} and Theorem~\ref{th1}, one
may be tempted to dream of a positive answer to the following question
raised by Falbel:\\

\noindent{\bf Question:} Let $M$ be a non-compact hyperbolic
3-manifold of finite volume, and let $\rho:\pi_1(M)\rightarrow
PU(2,1)$ be a \emph{discrete} representation such that every peripheral
subgroup is mapped to a cyclic group generated by a unipotent
element. Is the manifold at infinity of $\rho(\pi_1(M))$ homeomorphic
to $M$?\\

The requirement that peripheral $\Z\oplus\Z$ subgroups are mapped to
subgroups isomorphic to $\Z$ is included because of the results
in~\cite{falbelfigure8} (the corresponding representation is boundary
injective, but the domain of discontinuity is empty, so there is no
manifold at infinity at all).

The next result shows that the answer is negative in general; it
suggests one needs to be very cautious when studying
representations of $3$-manifolds into $PU(2,1)$.
\begin{thm}\label{thm:conjugacy}  
  The groups $\rho_{3\_3}(\pi_1(\verb|m015|))$ and
   $\rho_{4\_3}(\pi_1(\verb|m009|))$ are conjugate in $PU(2,1)$. 
\end{thm}
In section~\ref{sec:chford}, we will prove that both representations
are discrete, have non-empty domain of discontinuity, and that the
image group has only isolated fixed points in $\ch 2$.  In particular
at least one of the two manifolds gives a negative answer to Falbel's
question. Using Theorem~\ref{th1} (which we will not prove), one would see
that the negative answer is actually provided by the manifold
\verb|m015|. In other words, the manifold at infinity of
$\rho_{3\_3}(\pi_1(\verb|m015|))$ is \emph{the wrong manifold}, in the
sense that it is homeomorphic to \verb|m009|, \emph{not} to
\verb|m015|. 

Let us emphasize once more that Theorem~\ref{thm:conjugacy} shows
that a discrete, boundary unipotent representation of the fundamental
group of a given non-compact manifold into $PU(2,1)$, is not
necessarily the holonomy of a spherical CR uniformization of that
manifold, even if the peripheral holonomy is cyclic.

In section~\ref{sec:fillings} we will describe a technical feature
shared by all the noncompact hyperbolic manifolds that are known to
admit a spherical CR uniformization. For the time being, this feature
may serve as an explanation for the existence of these exotic
uniformizations. Specifically, it turns out that all finite volume non
compact 3-manifolds that are known to admit a spherical CR
uniformization all admit an exceptional Dehn filling that is a Seifert
fibration over a $p,q,r$-orbifold (see Theorem~\ref{thm:exceptional})
with $p,q,r\geq 3$. For 3-manifolds that do not have such Seifert Dehn
fillings, no satisfactory evidence of non-existence of spherical CR
uniformizations is presently available (but the FKR census gives no
discrete representation with cyclic unipotent boundary holonomy).

Now theorem~\ref{thm:conjugacy} will be an obvious consequence of a
stronger statement, namely Theorem~\ref{thm:trianglegroups}, which
gives a way to reconstruct the image of the two relevant FKR census
representations directly by deforming triangle groups. The idea is to
take the obvious embedding of the 3,3,5-triangle group, obtained via
the injection $SO(2,1)\subset SU(2,1)$, where reflections in $H^2_\R$
are extended to complex reflections in $H^2_\C$. Note that this
representation is type-preserving, in the sense that elliptic
(resp. parabolic, loxodromic) elements are mapped to elliptic
(resp. parabolic, loxodromic) elements.

The index two subgroup of words of even length in the triangle group
has a manifold at infinity which is a Seifert fibration over the
3,3,5-orbifold (see chapter~4 of~\cite{richBook}). It is well known
that, modulo conjugation in $PU(2,1)$, the deformation space of this
representation of the 3,3,5-triangle group is an interval
(see~\cite{schwartzICM} for instance), and one expects that, at least
for small deformations, the group should remain discrete, and the
manifold at infinity should not change its homeomorphism type.

The idea is then to consider the first place in the deformation space
where the representation is no longer type-preserving. A conjectural
quantitative analysis of when this happens is stated
in~\cite{schwartzICM}; in the case at hand, the first change in types
should occur when the word $I_2I_3I_1I_3$ (which is loxodromic in the
original triangle group) becomes parabolic. We will call the
corresponding group the first 3,3,5-triangle group with an
\emph{accidental parabolic element}, even though the validity of this
description really relies on the validity of Schwartz's
conjectures. This group is often denoted $(3,3,5;\infty)$ in the
literature.

In this paper, we describe explicit representations of
$\pi_1(\verb|m009|)$ and $\pi_1(\verb|m015|)$ onto the
$(3,3,5;\infty)$ triangle group, and show that they map peripheral
subgroups to cyclic subgroups generated by the accidental parabolic
element.  The accidental parabolic element is a unipotent element,
i.e. it has 1 as its only eigenvalue (see
section~\ref{sec:fillings}). From this one can easily identify these
two representations as specific representations in the FKR census (see
section~\ref{sec:trianglegroups}).
 
We finish by noting that the result of
Theorem~\ref{thm:trianglegroups}, even though it gives a way to bypass
the use of the FKR census, was widely inspired by detailed inspection
of the representations in the census. Originally, the author had computed
Ford domains for the image of the census representations that looked
discrete, and noticed that two of these Ford domains were isometric.

In that sense, the Ford domain of a group (centered at the fixed a
suitably chosen parabolic element) gives a very efficient conjugacy
invariant of the group, just like it does in the real hyperbolic case
(in that case, the tiling dual to the tiling by the Ford domain
produces the so-called canonical decomposition,
see~\cite{epsteinpenner}).

\begin{flushleft}{\bf Acknowledgements:} 
  This work was partly supported by the ANR, through the grant SGT
  (``Structures G\'eom\'etriques Triangul\'ees''). The author also
  benefited from generous support from the GEAR network (NSF grants
  DMS 1107452, 1107263, 1107367), via funding of an long term visit at
  ICERM; the author warmly thanks ICERM for its hospitality,
   Elisha Falbel, Pierre-Vincent Koseleff and Fabrice
    Rouillier for sending him an early version of their census; Ben
    Burton, Nathan Dunfield and Dave Futer for useful conversations
    related to this work; and finally the referee, who suggested
    several improvements in the manuscript.
\end{flushleft}

\section{Ford domains in $H^3_\mathbb{R}$} \label{sec:realford}

We briefly recall the general notion of Ford domain for discrete
subgroups of $PSL_2(\C)=Isom(H^3_\R)$, and describe these domains for
the special case of the holonomy group of the complete hyperbolic
structures on three specific 3-manifolds, namely \verb|m004|,
\verb|m009|, \verb|m015| in the Hildebrand-Weeks census.

\subsection{Real hyperbolic space and Ford domains}

Here we view $H^3_\R$ as the upper half space 
$$
\{(z,t)\in \C\times \R: t>0\},
$$
with the metric $(dx^2+dy^2+dt^2)/t^2$.
Recall that $\Gamma\subset PSL_2(\C)$ acts on $\partial_\infty
H^3_\R\simeq P^1_\C$ by its linear action on $C^2$. Working in an
affine chart, one gets an action by fractional linear transformations,
i.e.
$$
z\mapsto \frac{az+b}{cz+d}. 
$$ 
The basic point is that these maps extend to give an isometric
action of $PSL_2(\C)$ on $H^3_\R$, uniquely determined by the fact
that circles in $\C$ give the boundary of a unique sphere in $\C\times
\R$ orthogonal to the horizontal plane $\C$.

A formula for the extension can be most easily obtained by seeing
$H^3_\R$ as a totally geodesic subspace of $H^4_\R$, which is
isometric to $H^1_\H$, the 1-dimensional quaternionic hyperbolic
space, see~\cite{parkerjyva} for instance.
Concretely, one gets
$$ 
(z,t)\mapsto \left(\frac{(az+b)(\bar c \bar z+\bar d)+a\bar c
  t^2}{|cz+d|^2+|c|^2t^2},\frac{t}{|cz+d|^2+|c|^2t^2}\right).
$$

Now let
$$
\gamma=\left(\begin{matrix}a&b\\c&d
\end{matrix}\right),
$$
and suppose $c\neq 0$, or equivalently that $\gamma$ does not fix
infinity.
\begin{dfn}
  The \emph{isometric circle} of $\gamma$ is the circle in $\C$ where
  the derivative of the corresponding fractional linear transformation
  has modulus 1, which has equation $|cz+d|=1$. The \emph{isometric
    sphere} of $\gamma$ is the unique sphere in $\C\times \R$
  orthogonal to $\C$ that contains its isometric circle, with equation
  $|cz+d|^2+|tc|^2=1$.
\end{dfn}
Note that the circle has center $-d/c$ and radius $1/|c|$. Moreover,
it is easy to see that $\gamma$ maps its isometric circle to the
isometric circle of $\gamma^{-1}$. Finally, note that the hemispheres
in $\C\times\R$ orthogonal to $\C$ are totally geodesic copies of
$H^2_\R$ in $H^3_\R$.
\begin{dfn}
  The Ford domain of $\Gamma$ in $\H^3_\R$ is the connected component
  containing infinity of the complement of all isometric spheres of
  elements in $\Gamma$.
\end{dfn}
Provided the group $\Gamma$ is discrete, its Ford domain is a
fundamental domain for the action of $\Gamma$ if and only if no
element of $\Gamma$ fixed infinity.
It is useful to normalize the matrices so that infinity does have a
nontrivial stabilizer, in which case the stabilizer acts on $\C\simeq
\partial_\infty H^3_\R$ by a group $S$ of similarities, and one gets a
fundamental domain for the decomposition of $\Gamma$ into $S$-cosets
(see~\cite{beardon} for instance).

\subsection{The Hildebrand-Weeks census}

The Hildebrand-Weeks census is a list of all 1-cusped hyperbolic
3-manifolds that can be built by gluing up to 5 tetrahedra,
see~\cite{hildebrand-weeks}. For completeness, we mention that the
census was subsequently extended to allow for more tetrahedra,
see~\cite{callahan-hildebrand-weeks}, see also the work of
Thistlethwaite~\cite{thistlethwaite}, and Burton~\cite{burton}, but
the manifolds we consider in this paper require only 3 ideal
tetrahedra.

Recall that when it exists, the complete hyperbolic metric actually
has finite volume, so the metric is unique by Mostow rigidity; in
other words, the Kleinian groups are determined \emph{up to
  conjugation} in $PSL_2(\C)$. Because of the issue of possible
conjugation, it is sometimes difficult to compare different groups,
but there is a canonical way to associate a triangulation, obtained by
taking the decomposition \emph{dual} to the Ford domain
(see~\cite{epsteinpenner}).

This canonical decomposition is of course encoded in SnapPy; we will
start our description from the output of SnapPy for each of the three
manifolds considered in this paper, obtained with the \verb|canonize|
command (we are using SnapPy version 2.0, but the commands we use are
so standard that they should remain stable in subsequent versions).
In order to avoid cumbersome notation, throughout
section~\ref{sec:realford}, we will use the same notation for
generators in group presentations, and their image in $PSL_2(\C)$, so
$x$ will sometimes stand for $\rho(x)$; this is reasonable because all
the representations we consider in this section are known to be faithful.

\subsubsection{\bf m004}

{\bf Presentation:}
$$
\langle\ x,y\ |\ x[x,y][y^{-1},x^{-1}]\ \rangle
$$
{\bf Generators of a peripheral subgroup:}
$$
xy,\quad [x,y^{-1}][x^{-1},y^{-1}]
$$
{\bf Shape of the cusp:} 
$$
2i\sqrt{3}
$$
{\bf Triangular generators:}
$$
x^2yx^{-1} = \left(\begin{matrix} 1 & 0 \\ \alpha & 1\end{matrix}\right),\quad 
xyx = \left(\begin{matrix} 1 & \alpha \\ 0 & 1\end{matrix}\right),
$$ 
where $\alpha=\frac{i-\sqrt{3}}{2}$, which has minimal polynomial
  $x^4-x^2+1$. Note that there is in fact a representation into $PSL_2(K)$ for a
smaller number field, but we will not need this here.\\

The triangular shape of $s=x^2yx^{-1}$ and $t=xyx$ is easy to guess
from the SnapPy canonical presentation; it is not completely obvious
that these two matrices generate the whole group, so we mention that
$$
[s^{-1},t]
=x[y^{-1},x^{-1}]x^2yx^{-2}y^{-1}x^{-1}=x[y,x]xyx^{-2}y^{-1}x^{-1} = xyt^{-1},
$$
so $xy$ is in the group generated by $s$ and $t$, which easily proves
that $s$ and $t$ generate.

Finally, we summarize how to obtain the minimal polynomial of
$\alpha$. Using $s$ and $t$ as generators, the previous discussion
implies
$$
x=\left(\begin{matrix} -\alpha^2+1 & -\alpha^3 \\ \alpha^3 & \alpha^4+\alpha^2+1
\end{matrix}\right),\quad 
y=\left(\begin{matrix} z\alpha^4+2\alpha^2+1 & 2\alpha^3+\alpha \\ -2\alpha^3 & -2\alpha^2+1
\end{matrix}\right).
$$ 
The relation in the presentation translates into a set of
polynomial equations in a. Specifically, we require that
$M = x[x,y][y^{-1},x^{-1}]$
is scalar, which becomes
$$
\begin{array}{c}
M_{1,2}
= -\alpha^3(\alpha^4+\alpha^2+1)(\alpha^4-\alpha^2+1)(\alpha^8-\alpha^4+2\alpha^2+1) = 0\\
M_{2,1}
= \alpha^3(\alpha^4+\alpha^2+1)^2(\alpha^4-\alpha^2+1)^2 = 0\\
M_{1,1}-M_{2,2}
= -\alpha^2(\alpha^4+\alpha^2+1)(\alpha^4-\alpha^2+1)(\alpha^{10}+2\alpha^8-\alpha^6+2\alpha^4+\alpha^2+2)
\end{array},
$$ 
so taking $\alpha$ to be a root of $x^4-x^2+1$ will certainly give
a representation into $PSL_2(\C)$ (other choices of $\alpha$ will give
Galois conjugate representations).

In particular, for $\alpha=(i-\sqrt{3})/2$, one gets lower triangular
matrices for the stabilizer of a cusp by computing
$$
x^2yx^{-1} = \left(\begin{matrix} 1 & 0 \\ \alpha & 1 
\end{matrix}\right),\quad 
x[x,y][x,y^{-1}]x^{-1} = \left( \begin{matrix} -1 & 0 \\ -\sqrt{3}-3i & -1 
\end{matrix}\right).
$$ 
and the ratio of the lower left entries give the shape of the cusp, namely
$$
-\frac{\sqrt{3}+3i}{\alpha} = 2i\sqrt{3}.
$$

\subsubsection{\bf m009}

{\bf Presentation:}
$$
\langle\ x,y\ |\ x[x,y]x[x,y^{-1}]\ \rangle
$$
{\bf Generators of a peripheral subgroup:} 
$$
xy,\quad x^{-1}y^{-1}x^3y^{-1}x^{-1}y
$$
{\bf Shape of the cusp:} 
$$
i\sqrt{7}
$$

The SnapPy representations give matrices with entries in
$\mathbb{Q}(i,\beta)$, where
$$
\beta^4+\beta^2+2=0.
$$ 
The specific matrices are then given by
$$
x=\left(\begin{matrix} -\beta^3-\beta & i \\
                             -i             & \beta
              \end{matrix}\right),\quad 
y=\left(\begin{matrix} -\beta^3 & i \\
                             -i(\beta^2+1)  & \beta
              \end{matrix}\right),
$$ 
where $\beta$ is the root which is given to eight decimal places by
$\beta_0=0.67609672 + 0.97831834i$. Note that the other roots of the
polynomial $\beta^4+\beta^2+2$, namely $-\beta_0$ and $\pm
\overline{\beta}_0$, happen to give representations that are conjugate
in $PSL_2(\C)$ either to the above representation or to its complex
conjugate, even though this is far from a general phenomenon (Galois
conjugates of lattice representations are often non-discrete).

It is easy to see that the single fixed point of $xy$ is
$-i\beta$, so that the matrix
$$
q = \left(\begin{matrix}
  -i\beta & 0\\ 1&1/\beta^2
\end{matrix}\right)
$$
conjugates $xy$ into 
$$
\left(\begin{matrix}1&1\\0&1
\end{matrix}\right),
$$
and then the other generator of the above peripheral subgroup gets conjugated to
$$
\left(\begin{matrix}
  -1 & 1+2\beta^2\\
  0  & -1
\end{matrix}\right).
$$
Note that
$$
(2\beta^2+1)^2=4\beta^4+4\beta^2+1=-7,
$$ 
so that $1+2\beta^2=\pm i\sqrt{7}$, and with the above choice of
the root $\beta$ given above, one checks it is actually $i\sqrt{7}$.
In any case, the cusp section corresponds to a square lattice,
generated by $1$ and $i\sqrt{7}$.

\subsubsection{\bf m015}

{\bf Presentation:}
$$
\langle\ x,y\ |\ [x,y^{-1}]x^3[y,x^{-1}]y^2\ \rangle
$$
{\bf Generators of a peripheral subgroup:}
$$
xy,\quad (xy)^2[x,y^{-1}]x[y^{-1},x]y^{-1}xy
$$
{\bf Shape of the cusp:} 
$$
4(\gamma-1)
$$
{\bf Triangular generators:}
$$
yx = \left(\begin{matrix} -1 & -\gamma \\ 0 & -1\end{matrix}\right),\quad 
xyx = \left(\begin{matrix} 1 & 0 \\ \gamma & 1\end{matrix}\right),
$$ 
where $\gamma$ is a complex root of $x^3-x^2+1$ (say the one which
  is approximately $0.87743883-0.74486176i$).\\

The shape of the matrices for the triangular matrices $yx$ and $xyx$
is easily guessed from the SnapPy matrices (which gives only numerical
approximations of those matrices), and the minimal polynomial for
$\gamma$ is a consequence of the relation given for these
generators. Note also that it is clear that $xy$ and $xyx$ generate
the group.

\subsection{Ford domains in $H^3_\R$}

The above information allows to study the Ford domains in
$H^3_\R$. The rough idea is to compute a symmmetric set of somewhat
short words in the generators, to consider their isometric spheres (as
well as their images under the cusp group), which gives a ``partial''
Ford domain (i.e. a polytope that may or may not be equal to the
actual Ford domain).

In order to check whether the partial Ford domain is equal to the Ford
domain, one can then apply the Poincare polyhedron theorem in order to
check whether copies of the Ford domain under the group tile
$H^3_\R$. The fact that our set of words is symmetric implies that
faces of the Ford domain are paired (the face for $\gamma$ is sent to
the one for $\gamma^{-1}$), and one needs to check the cycle
conditions (which roughly say that on the level of codimension two
faces, the side-pairings induce a local tiling of $H^3_R$).

The Ford domains are by construction invariant under the action of the
peripheral subgroups, so they are not fundamental domains (they are
only fundamental domains for the decomposition of the group into
cosets of a given peripheral subgroup).

The classical way to obtain an actual fundamental domain for the
action of the group is to intersect it with a fundamental domain for
the action of the cusp group, which is a lattice in $\C$. Hence one
can simply intersect the Ford domain with a vertical prism over a
parallelogram in $\C$.

Another way is to select a representative for each orbits of faces of
the Ford domain under the action of the cusp group (we may assume that
the union of the face representatives is connected). A fundamental
domain is then obtained again as a vertical prism, but over a union
faces contained in spheres.  

\begin{figure}[htbp]
  \subfigure[m004]{\epsfig{figure=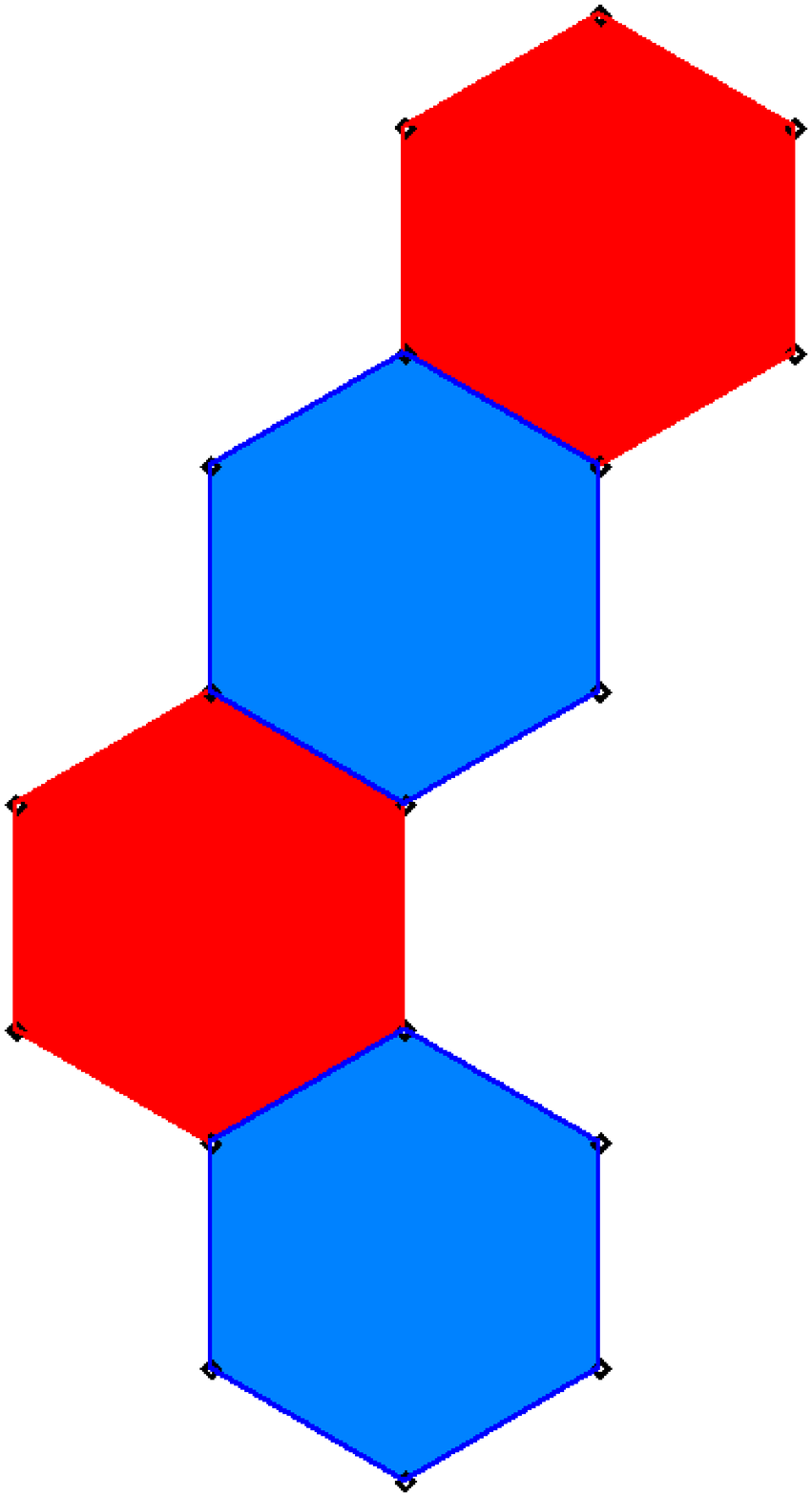,width=0.2\textwidth}}\quad
  \subfigure[m009]{\epsfig{figure=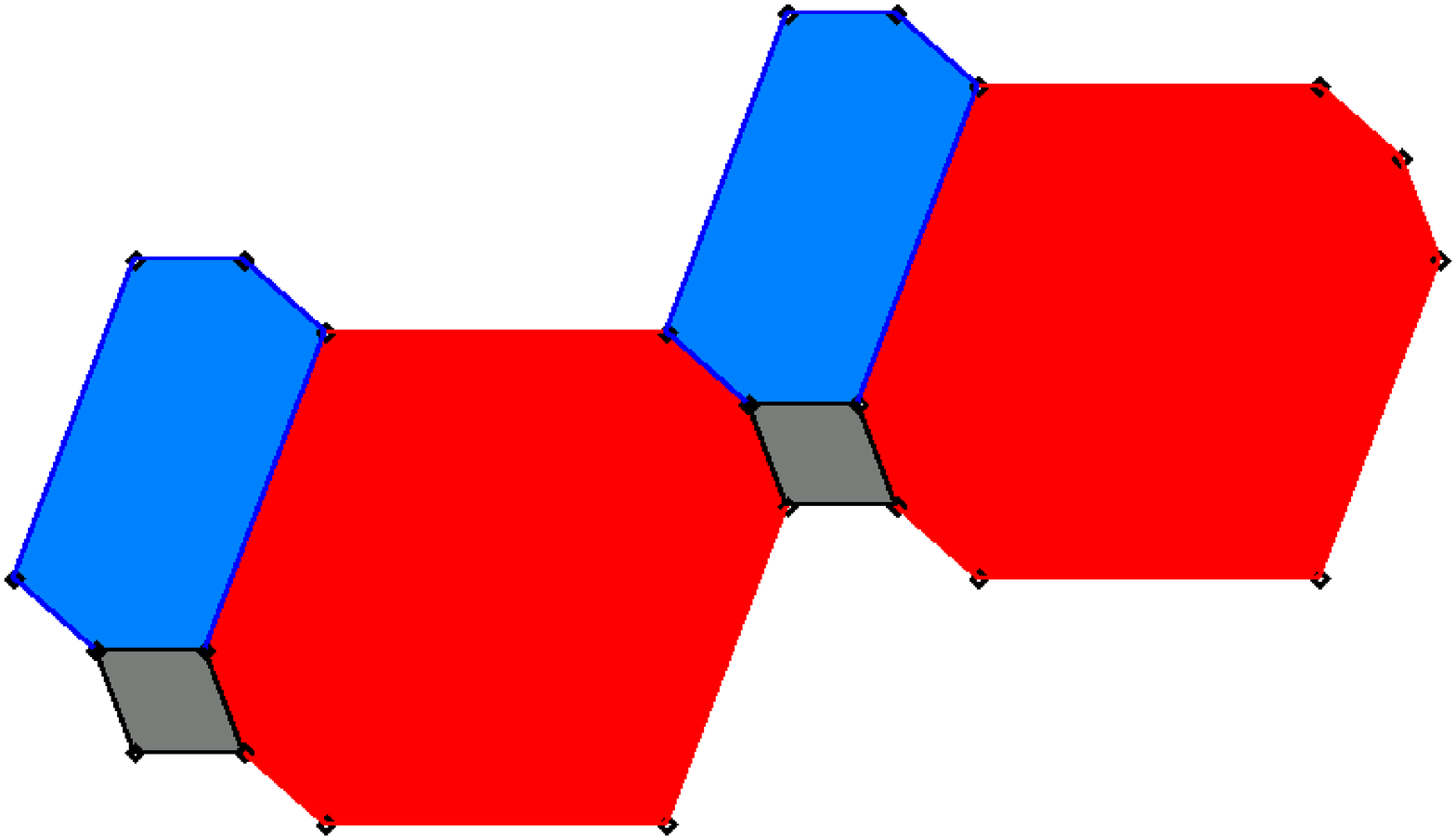,width=0.5\textwidth}}
  \subfigure[m015]{\epsfig{figure=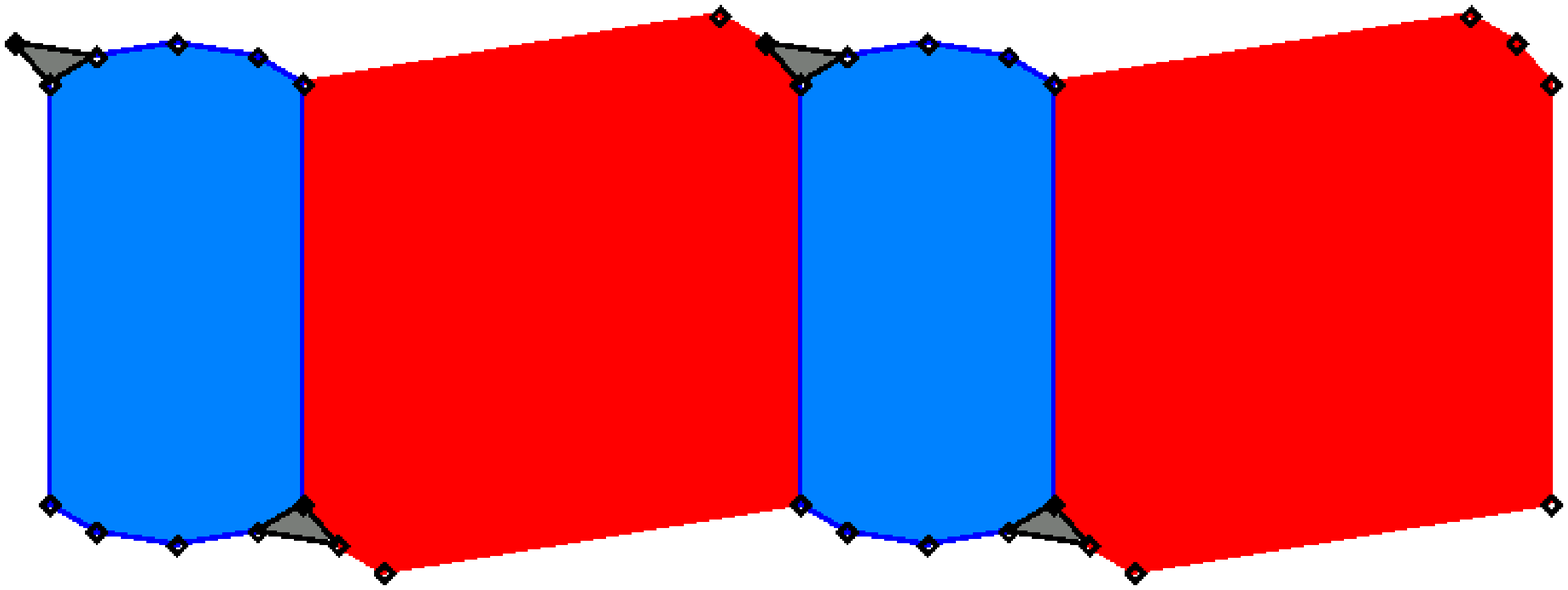,width=0.8\textwidth}}
  \caption{Prism description of some $1$-cusped hyperbolic
    manifolds}\label{fig:prisms}
\end{figure}

Looking at the picture from infinity, we see polygons in $\C$, which
are depicted in Figure~\ref{fig:prisms} for the three groups that
appear in this paper.

The picture for \verb|m004| is very classical, and appears already
in~\cite{rileyFigure8}. In that case all isometric spheres bounding
the Ford domain have the same radius, and they can be taken to be
given simply by the spheres of radius 1 centered at points in the
Eisenstein lattice $\Z[\omega]$, $\omega=(-1+i\sqrt{3})/2$.

In the other two cases, there are three different radii for the
isometric spheres that bound the Ford domain. The pictures can of
course be obtained directly in SnapPy (using the
\verb|cusp_neighborhood| command).

\section{Exceptional fillings of census manifolds}\label{sec:fillings}

It is well known that for every complete finite volume hyperbolic
manifold $M$, all but finitely many Dehn fillings are hyperbolic,
see~\cite{thurstonNotes},~\cite{neumannzagier},~\cite{petronioporti}. For
simplicity, we assume that $M$ only has one cusp, and write $M(p/q)$
for the Dehn filling of slope $p/q$.

Here $p/q$ is a reduced fraction (or possibly infinity), and is
supposed to specify the slope of a circle on the boundary torus that
is chosen to bound a meridian in the 2-torus that gets glued to $M$.
Note that $M(p/q)$ is only well-defined once a meridian and longitude
have been chosen in the boundary torus of $M$, which can be done
canonically when $M$ is a knot complement (otherwise we will use the
meridian and longitude provided by SnapPy when using the
Hildebrand-Weeks census with the notation \verb|m0jk| that we have
used throughout the paper).

A lot of work has been done over the years in order to make the
``finitely many'' part of the statement effective, i.e. to give bounds
on the number of all exceptional slopes for knot complements, or to
characterize possible Dehn fillings. For instance, exceptional
surgeries of 2-bridge knots are classified in~\cite{brittenhamwu}. A
convenient place to find a list for a lot of manifolds in the
beginning of the Hildebrand-Weeks census is~\cite{magicmanifold}. One
can also gather a lot of experimental evidence for their statements by
using recent versions of SnapPy (simply produce a triangulation for
various Dehn-filling, and pass them on to Regina for further
analysis).

A key observation, in view of the main statement
of~\cite{derauxfalbel} (see also Table~\ref{tab:discrete}) is that
\verb|m004|, \verb|m009| and \verb|m015| all have exceptional Dehn
fillings that produce Seifert fibrations over spherical
$(p,q,r)$-orbifolds (these are often called \emph{small} Seifert
fibrations), all with $p,q,r\geq 3$. For convenience of notation, we
write $M_4$ (resp. $M_9$, $M_{15}$), for \verb|m004|, \verb|m009| and
\verb|m015|.
\begin{thm}\label{thm:exceptional}
  The following Dehn fillings are small Seifert spaces.
  \begin{itemize}
  \item $M_4(\pm 3)$ is a Seifert fibration over the $(3,3,4)$-spherical
  orbifold.
  \item $M_9(-2)$ is a Seifert fibration over the $(3,3,5)$-spherical
  orbifold.
  \item $M_{15}(1)$ is a Seifert fibration over the $(3,3,5)$-spherical
  orbifold. 
  \end{itemize}
\end{thm}

These claims can be gathered somewhat painfully from the information
in~\cite{magicmanifold}, see Proposition~2.2~(26) of that paper. The
main difficulty is that the slope of a Dehn filling depends on the
basis of the homology that is used. In this paper, we use the basis
that is in the SnapPy database, whereas Martelli-Petronio use the
basis induced from the canonical bases for the homology of the three
boundary components of the magic link complement. The relationship
between the two follows from elementary Kirby calculus
(see~\cite{rolfsen}, p.265).

Note that Table~8 of~\cite{magicmanifold} also contains the claim
about $M_4=N(1,2)$ and $M_9=N(1,3)$ ($N(p,q)$ denotes a Dehn filling
of two of the three cusps of the magic manifold with slopes $p$ and
$q$ respectively); $M_9$ is not a knot in $S^3$, but it is a knot in
$\mathbb{R}P^3$.

For the reader who is more versed in geometry than topology, the best
way to check Theorem~\ref{thm:exceptional} is to use SnapPy in
conjunction with Regina. Indeed, the following SnapPy commands will
produce a triangulation for the Dehn filling of \verb|m009| with slope
$-2$:

\begin{center}
  \begin{tabular}{l}
    \verb|M = Manifold('m009')|\\
    \verb|M.dehn_fill((-2,1))|\\
    \verb|T = M.filled_triangulation()|\\
    \verb|T.save()|
  \end{tabular}
\end{center}
This triangulation can then be imported in Regina; either the manifold
is recognized right away, or it can be recognized by performing a
Census Lookup (indeed for the Dehn fillings that appear in the present paper,
this lookup turns out to be successful). In the above example, Regina
describes the Dehn filling by
\begin{center}
  \verb| SFS [S2: (3,1) (3,1) (5,-4)]|
\end{center}
where SFS stands for Seifert Fibered Space, and the rest gives gluing
information. The base of the fibration is $S^2$, and there are three
singular fibers, with gluing data given by the following pairs of
integers. For the precise meaning of the gluing data, see chapter 12
of~\cite{hempel}. Here we simply mention that the base is a sphere
with three orbifold points, with weights 3,3 and 5. In particular, by
contracting all the fibers, one gets a homomorphism of $\pi_1(M_9)$
onto a $(3,3,5)$-triangle group.

Apart from the possibility of a bug in these well-established computer
programs, this computer check can be regarded as a proof because it is
purely combinatorial in nature.  Note also that we have stated
Theorem~\ref{thm:exceptional} only for motivational purposes. It
implies that the fundamental groups of $M_4$, $M_9$ and $M_{15}$ all
admit homomorphisms onto a $(3,3,n)$-triangle group, with $n=4$ or
$5$; in fact, in the next section, we will construct explicit such
homorphisms without appealing to Theorem~\ref{thm:exceptional}.

\section{Complex hyperbolic geometry and triangle groups}\label{sec:trianglegroups}

The main goal of this section is to give a triangle group
interpretation of two of the discrete groups that occur as holonomy
groups in the FKR census, namely $\rho_{4\_3}(\pi_1(M_9))$ and
$\rho_{3\_3}(\pi_1(M_{15}))$, see
Theorem~\ref{thm:trianglegroups}. This identification will immediately
imply that these holonomy groups are actually conjugate to each other
in $PU(2,1)$.

For basics on complex hyperbolic geometry and triangle groups,
see~\cite{goldman},~\cite{schwartzICM},~\cite{derauxfalbel} for instance.
Recall that complex hyperbolic triangle groups generated by three
complex involutions $I_1$, $I_2$, $I_3$ that satisfy
$$
(I_1I_2)^p=(I_2I_3)^q=(I_3I_1)^r=Id.
$$ 
In that context, the condition $p,q,r\geq 3$
can be thought of as requiring that the triangle does not have any
right angle, or equivalently that the corresponding $\R$-Fuchsian
triangle group admit non-trivial deformations:
\begin{prop}
  $(2,q,r)$-triangle groups in $PU(2,1)$ are rigid, but
  $(p,q,r)$-triangle groups with $p,q,r\geq 3$ have a 1-dimensional
  character variety.
\end{prop}

\begin{pf}
  This follows from the explicit parametrization of triangle
  groups. Given an irreducible triangle (i.e. without any global fixed
  point), we can take three polar vectors $v_1,v_2,v_3$ to the mirrors
  of generating involutions as a basis of $\C^3$. After suitably
  normalizing these vectors, we may assume $\langle v_j,v_j\rangle=1$
  for all $j$, and we may also assume $\langle v_1,v_2\rangle$ and
  $\langle v_2,v_3\rangle$ are real (but in general, $\langle
  v_1,v_2\rangle$ will not be real). An invariant of the phase change
  for the $v_j$'s is given by the triple Hermitian product
  $$
  \langle v_1,v_2\rangle \langle v_2,v_3\rangle \langle v_3,v_1\rangle,
  $$ 
  whose argument is sometimes called the angular invariant of the
  triangle; one checks that for every $(p,q,r)$, only an interval of
  values of the angular invariant can be realized by complex
  hyperbolic triangles, characterized by requiring that the Hermitian matrix
  \begin{equation}\label{eq:nondiagonal}
  H=\left(\begin{matrix}
    1 & -\cos\frac{\pi}{r} & -\cos\frac{\pi}{q}\varphi\\
    -\cos\frac{\pi}{r} & 1 & -\cos\frac{\pi}{p}\\
   -\cos\frac{\pi}{q}\overline{\varphi} & -\cos\frac{\pi}{p} & 1
  \end{matrix}\right),
  \end{equation}
  have negative determinant, where $\varphi=e^{it}$ and $t\in\R$ is the
  angular invariant. Note that $t=0$ corresponds to $\R$-Fuchsian groups.

  When none of the cosines is zero, the signature condition translates
  to a lower bound on $\cos t$, which gives an interval containing 0
  of admissible values of $t$. If one of the cosines is 0, then all
  $(p,q,r)$-triangle groups are $\R$-Fuchsian, and they are all
  isometric to each other.
\end{pf}

The following result is a strengthening of the claim that the discrete
representations of $M_9$ and $M_{15}$ in the FKR census are
conjugate. The idea is to identify the image of these representations as explicit
triangle groups; the fact that these two 3-manifolds both admit an
exceptional Dehn filling that is a Seifert fibration over a
3,3,5-orbifold (see Theorem~\ref{thm:exceptional}) immediately implies that
their fundamental group admits a homomorphism onto the 3,3,5-triangle
group (obtained by contracting of the fibers).

In Theorem~\ref{thm:trianglegroups}, we describe an explicit such
homomorphism and show that the corresponding peripheral holonomy is
cyclic unipotent; hence the corresponding representations must
actually appear somewhere in the FKR census, and they are easily
identified in the census list using the field of cross ratios.

\begin{thm} \label{thm:trianglegroups}
  \begin{enumerate}
  \item Up to conjugacy and complex conjugation, there is a unique
    $(3,3,4)$-triangle group such that $I_2I_3I_1I_3$ is
    parabolic. Its even length words subgroup is conjugate to both
    $\rho_{1\_1}(\pi_1(M_4))$ and $\rho_{4\_1}(\pi_1(M_4))$, or in
    other words to the holonomy group of the unique boundary unipotent
    spherical CR uniformization of the figure eight knot complement.
  \item Up to conjugacy and complex conjugation, there is a unique
    $(3,3,5)$-triangle group such that $I_2I_3I_1I_3$ is
    parabolic. The even length subgroup of that triangle group is
    conjugate to both $\rho_{4\_3}(\pi_1(M_9))$ and
    $\rho_{3\_3}(\pi_1(M_{15}))$.
  \end{enumerate}
\end{thm}

\begin{rmk}
  The triangle groups that appear in Theorem~\ref{thm:trianglegroups}
  are often denoted $(3,3,4;\infty)$ and $(3,3,5;\infty)$,
  respectively.
\end{rmk}

\begin{pf}
  We treat the case of $(3,3,5)$-triangles, the other one being 
  entirely similar. In that case, the matrix of equation~\eqref{eq:nondiagonal} reads
\begin{equation}\label{eq:hermform}
  H = \left(\begin{matrix}
    1 & -\frac{1}{2} & -\frac{1+\sqrt{5}}{4}\varphi\\
    -\frac{1}{2} & 1 & -\frac{1}{2}\\
   -\frac{1+\sqrt{5}}{4}\overline{\varphi} & -\frac{1}{2} & 1
  \end{matrix}\right),
\end{equation}
  which has negative determinant for $\varphi=e^{it}$ and
  $|t|<\arccos\frac{\sqrt{5}-3}{2}$.
  
  It is not difficult to verify that $1$ is always an eigenvalue of
  $I_2I_3I_1I_3$, and that it has real trace (the latter
    follows from the fact that it is the product of two involutions,
    namely $I_2$ and $I_3I_1I_3$, so it is conjugate to its own
    inverse).  In particular, if it is parabolic, then it must
  actually be unipotent, hence its trace must be equal to 3.

  Now even when $I_2I_3I_1I_3$ is not parabolic, its two other
  eigenvalues are complex conjugate (provided we work with matrices in
  $SU(2,1)$), and their sum is
  $$
  \frac{1+\sqrt{5}}{2}(2c+1),
  $$ 
  which is equal to 2 for $c=\sqrt{5}/2-1$. This corresponds to taking 
  $$
  \varphi = ( \sqrt{5}-2+i\sqrt{4\sqrt{5}-5} )/2.
  $$ which is one of the complex roots of the polynomial
  $x^4+4x^3+x^2+4x+1$. In other words, the relevant triangle group can
  be written in terms of matrices with entries in
  $K=\Q(\varphi)$\label{pagenbfield}, which is number field of
  degree 4 (beware that this extension is not Galois). The matrices
  actually have entries in the ring of integers $\mathcal{O}_K$
  (recall that generators are complex reflections of order 2).

  Note that the above value of $\varphi$ is indeed in the range where
  the Hermitian form has signature $(2,1)$. In fact, for that
  value of $c$, one gets
  $$
  \det(H) = - \frac{1+\sqrt{5}}{16}.
  $$ 

  \noindent
  {\bf The case of \verb|m009|}

  From SnapPy, we gather that $\pi_1(M_9)$ has the presentation
  $$ 
  \langle\ a,b,c,d\ |\ bac=db, c=ad, ca^{-1}bd^{-1}=id\ \rangle,
  $$
  with a peripheral subgroup generated by $b^{-1}adc^{-1}d$ and
  $d^{-1}cd^{-1}bc^{-1}db^{-1}$. This can be simplified to
  $$
  \langle a,d\ |\ a^2[a,d][a,d^{-1}]\ \rangle,
  $$ 
  with a peripheral subgroup generated by $[d^{-1},a]d$ and
  $d^{-1}a[a,d^{-1}]a^{-1}$.

  We describe an explicit homomorphism from $\pi_1(M_9)$ to the
  $(3,3,5)$ triangle group, which maps into the index two subgroup of
  even length words.
  $$
  \langle\ I_1,I_2,I_3\ |\ (I_1I_2)^3, (I_2I_3)^3, (I_3I_1)^5\ \rangle.
  $$ 
  Hoping that no confusion will arise, we use word notation, so that
  we write $123212$ for $I_1I_2I_3I_2I_1I_2$, for instance.

The map
  $$
  \left\{\begin{array}{l}
  a\mapsto 2132\\
%  b\mapsto 132132\\
%  c\mapsto 21321232\\
  d\mapsto 1232
  \end{array}\right.
  $$ 
extends to a homomorphism $\sigma:\pi_1(M_9)\rightarrow
\Gamma(3,3,5;\infty)$.  We skip the routine verification of this
statement. Note that, under this homomorphism, it is routine as well
to verify that the peripheral subgroups get mapped to cyclic groups
generated by a single unipotent element. 
%Specifically, one checks that
%$b^{-1}adc^{-1}d$ gets mapped to $2313$, whereas
%$d^{-1}cd^{-1}bc^{-1}db^{-1}$ gets mapped to $(2313)^2$, 
Specifically, one checks that $[d^{-1},a]d$ gets mapped to $2313$,
whereas $d^{-1}a[a,d^{-1}]a^{-1}$ gets mapped to $(2313)^2$, so that
our representation has cyclic unipotent boundary holonomy.\\

\noindent
{\bf The case of \verb|m015|}

We now sketch the corresponding arguments for
  $\pi_1(M_{15})$. The geometric presentation from SnapPy has the form
$$
\langle\ a,b,c,d\ |\ bad,\ cb^{-1}abd^{-1},\ cdc^{-1}a\ \rangle,
$$
which can easily be simplified to
$$
\langle\ a,b\ |\ b=ab^2a^{-1}[b^{-2},a^{-1}]\ \rangle.
$$
SnapPy also gives generators for a peripheral subgroup, namely
$d^{-1}c = b^{-1}a^{-1}b$ and $b^{-1}acbd^{-1}a^{-1}d = b^{-3} a^{-1} b^3
a^{-1} b^{-1}$.

One checks that 
$$
\left\{\begin{array}{l}
a\mapsto 2313\\
b\mapsto 1313\\
\end{array}\right.
$$ 
induces a well-defined homomorphism, and it maps the above two peripheral
elements to $131(2313)131$ and its inverse, respectively.\\

Now the homomorphisms we just constructed are both boundary
unipotent, these representations must appear somewhere in the FKR
census.  Up to complex conjugation, there are three representations of
$\pi_1(M_9)$ in the census. The field generated by the cross ratios of
the corresponding tetrahedra are different, only one has degree 4,
namely that for $\rho_{4\_3}$. Similarly, there is only one representation of $\pi_1(M_{15})$ with
the same field of cross ratios, namely $\rho_{3\_3}$.
\end{pf}

\section{The complex hyperbolic Ford domain}\label{sec:chford}

The Ford domain for boundary unipotent spherical CR uniformization of
the figure eight knot complement is studied
in~\cite{derauxdeformfig8}. We now sketch the corresponding analog for
the manifold \verb|m009| (giving all details would be much longer than
for the figure eight knot complement, but in fact it is not more
difficult).

We quickly review some basic material about Ford domains. It is
convenient to work with coordinates where the Hermitian form is given
by
\begin{equation}\label{eq:stdform}
  J=\left(\begin{matrix}
    0&0&1\\
    0&1&0\\
    1&0&0
  \end{matrix}\right),
\end{equation}
and we write $p_\infty$ for $(0,0,1)$. Recall that $\partial_\infty
H^2_\C\setminus \{p_\infty\}$ identifies with the Heisenberg group
$\C\times\R$, see section~4.2.2 of~\cite{goldman}. This identification
is obtained by considering the unipotent stabilizer of $p_\infty$;
note that a lower triangular matrix with unit diagonal is an isometry
of $J$ if and only if it has the form
$$
M(z,t)=\left(
\begin{matrix}
  1&0&0\\
  z&1&0\\
  -|z|^2/2+it&-\overline{z}&1
\end{matrix}
\right),
$$ 
for some $z\in\C$ and $t\in \mathbb{R}$. The group of such matrices acts
simply transitively on $\partial_\infty H^2_\C\setminus \{p_\infty\}$,
which gives the identification.
One easily computes that
$$
M(z,t)M(w,u)=M(z+w,t+u+Im(z\overline{w})),
$$ 
which suggests a group law on $\C\times\R$; up to a coefficient of 2,
this is the Heisenberg group law used in~\cite{goldman}, so we call
$(z,t)$ Heisenberg coordinates.

For a subgroup $\Gamma\subset PU(2,1)$, the Ford domain
$F_{\Gamma,p_\infty}$ centered at $p_\infty$ is given in homogeneous
coordinates by the set of vectors $Z\in \C^3$ that satisfy
$$
|\langle P_\infty,Z\rangle| \leq |\langle \tilde{g}P_\infty,Z\rangle|
$$
for all $g\in\Gamma$, where $\tilde{g}$ denotes any matrix
representative of $g\in\Gamma$. For each $g\in\Gamma$ not fixing
$p_\infty$, the set of points satisfying
\begin{equation}\label{eq:defface}
  |\langle P_\infty,Z\rangle| = |\langle \tilde{g}P_\infty,Z\rangle|
\end{equation}
is a so-called \emph{bisector}, which we denote by $\mathcal{B}_g$. It
is a basic fact that the trace at infinity of any such bisector, seen
in Heisenberg coordinates, is a \emph{bounded} topological sphere,
called a \emph{spinal sphere}.

\begin{dfn}\label{dfn:face}
  The 3-dimensional polyhedron given by the intersection of the Ford
  domain $F_{\Gamma,p_\infty}$ with $\mathcal{B}_g$ will be denoted by
  $b_g$.
\end{dfn}

When $\Gamma$ is discrete and $p_\infty$ is not fixed by any
non-trivial element in $\Gamma$, the Ford domain is actually a
fundamental domain for the action of $\Gamma$. If $p_\infty$ has a
discrete stabilizer $P\subset\Gamma$, then the Ford domain is only a
fundamental domain for the decomposition of $\Gamma$ into $P$-cosets
(see~\cite{derauxparkerpaupert} or~\cite{derauxfalbel} for instance).
  
It is usually difficult to determine this set explicitly, even though
it is fairly accessible to experimentation. Indeed, the boundary of
this domain is made up of pieces of bisectors, and bisector
intersections are now fairly well understood, see~\cite{goldman} for
instance. The basic point is that pairs of bisectors that occur in a
Ford domain have connected intersection, diffeomorphic to a smooth
disk. This imporant fact is stated in Theorem~9.2.6 of~\cite{goldman}
(in Goldman's language, pairs of bisectors that contain faces of a
Ford domain are called \emph{covertical} bisectors).
 
If $\Gamma$ can be represented by matrices in a given number field (of
reasonably small degree), there are computational tools to certify the
combinatorics of 2-faces of the Ford domain (see~\cite{derauxdeformfig8}).

We now apply these general notions to the groups that are the images
of the relevant FKR representations. Recall that we gave a detailed
description of these representations in
section~\ref{sec:trianglegroups}; our description relied on a specific
Hermitian form, see equation~\eqref{eq:nondiagonal}. In those
coordinates, one can easily work out formulas for the matrices $I_1$,
$I_2$, $I_3$, namely
$$
I_1=\left(\begin{matrix}
  1 & -1 & \alpha\\
  0 & -1 & 0\\
  0 & 0 & -1
\end{matrix}\right),
\quad
I_2=\left(\begin{matrix}
  -1 & 0 & 0\\
  -1 & 1 & -1\\
   0 & 0 & -1
\end{matrix}\right),\quad
I_3=\left(\begin{matrix}
  -1 & 0 & 0\\
  0 & -1 & 0\\
  \overline{\alpha} & -1 & 1
\end{matrix}\right),
$$
where $\alpha=-1+(1+\sqrt{5})(1-is)/4$ and $s=\sqrt{4\sqrt{5}-5}$.

We start by conjugating these three matrices so that they preserve the
standard Hermitian $J$, see equation~\eqref{eq:stdform}, and so that
$I_2I_3I_1I_3$ becomes (lower) triangular.
This is done by an easy
linear algebra computation, we give one explicit possible conjugation,
namely
$$
Q = \left(\begin{matrix}
  \sqrt{2}+(3-\sqrt{5})(-5+is)/4\sqrt{2} & 0          & -\sqrt{2}+(1+\sqrt{5})(1-is)/4\sqrt{2}\\
  (2+(2-\sqrt{5})(-3+is))/4\sqrt{2} & \sqrt{2}/2 & -\sqrt{2}\\
  (1-\sqrt{5})(-1+is)/4\sqrt{2}  & 0          & -\sqrt{2}+(1+\sqrt{5})(-1+is))/4\sqrt{2}\\
\end{matrix}\right).
$$
One easily checks that $Q^*HQ=J$. Writing $\tilde{I}_k=Q^{-1}I_kQ$, one checks that
$$
\tilde{I}_2=\left(\begin{matrix}
                 -1 & 0 & 0\\
                 0 & 1 & 0\\
                 0 & 0 & -1
\end{matrix}\right),\quad 
\tilde{I}_2\tilde{I}_3\tilde{I}_1\tilde{I}_3=
\left(
\begin{matrix}
  1 & 0 & 0\\
  1 & 1 & 0\\
  -\frac{1}{2}&-1&1
\end{matrix}
\right).
$$ 
The matrices for $\tilde{I}_1$ and $\tilde{I}_3$ are much more
complicated, so we do not write them out. In what follows, we simply will
write $I_k$ for $\tilde{I}_k$, since no confusion should arise.

%From this point on, we focus on a specific group $\Gamma$, see
%Definition~\ref{dfn:group}, and we take the center of the Ford domain
%to be the fixed point of $a=2313$ (recall that this stands for
%$I_2I_3I_1I_3$).
\begin{dfn}\label{dfn:group}
Let $\Gamma$ be the even length subgroup of the
$(3,3,5;\infty)$-triangle group, i.e. the subgroup generated by
$I_1I_2$ and $I_2I_3$. We define $a=I_2I_3I_1I_3$, and write $A$ for
$a^{-1}$.
\end{dfn}
Recall that the fixed point of $a$ is $p_\infty=(0,0,1)$. We write
$F=F_\Gamma=F_{\Gamma,p_\infty}$ for the corresponding Ford domain.
By construction, it cannot be a fundamental domain for the action of
$\Gamma$, since it is invariant under the action of the cyclic group
generated by $a$. It has infinitely many faces, but there are only ten
$a$-orbits of faces. Representatives of these orbits are given by the
faces $b_g\subset\B_g$ for elements $g$ in Table~\ref{tab:corefaces}.
\begin{table}
$$
  32, 23; 2321, 1232; 12, 21; 232131, 131232; 32131232, 23213123.
$$
\caption{The list of group elements whose orbit points define ten core
  faces, i.e. representatives of each orbit of faces under the
  $a$-action.}\label{tab:corefaces}
\end{table}
Note that we do not give the shortest possible word in
Table~\ref{tab:corefaces}, because we consider only even length words
in the triangle group. For instance, $32(e_3)=3(e_3)$,
$12(e_3)=1(e_3)$, etc.

We will number the elements in Table~\ref{tab:corefaces} as
$g_1=32,g_2=23,g_3=2321,\dots,g_{10}=23213123$ (in the same order as
listed in the table), and number the images of these under powers of
$a$ by setting for $j=1,\dots,10$ and $k\in\N$,
\begin{equation}\label{eq:numbering}
  a^k g_j = g_{10(2k-1)+j};\quad a^{-k} g_j = g_{20k+j}.
\end{equation}

\begin{figure}[htbp]
\hfill\subfigure[32]{\epsfig{figure=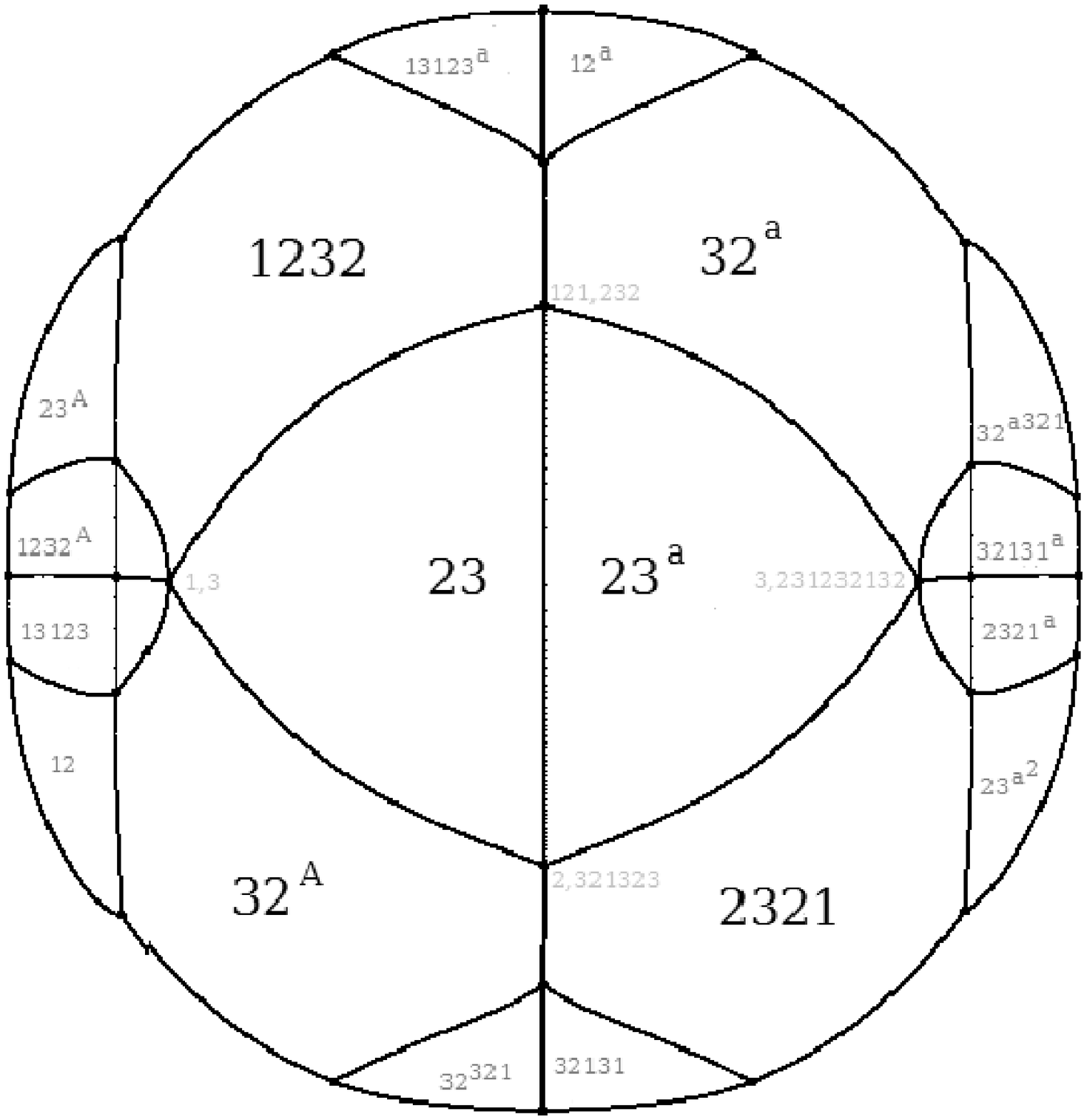, width=0.45\textwidth}}\hfill
\subfigure[23]{\epsfig{figure=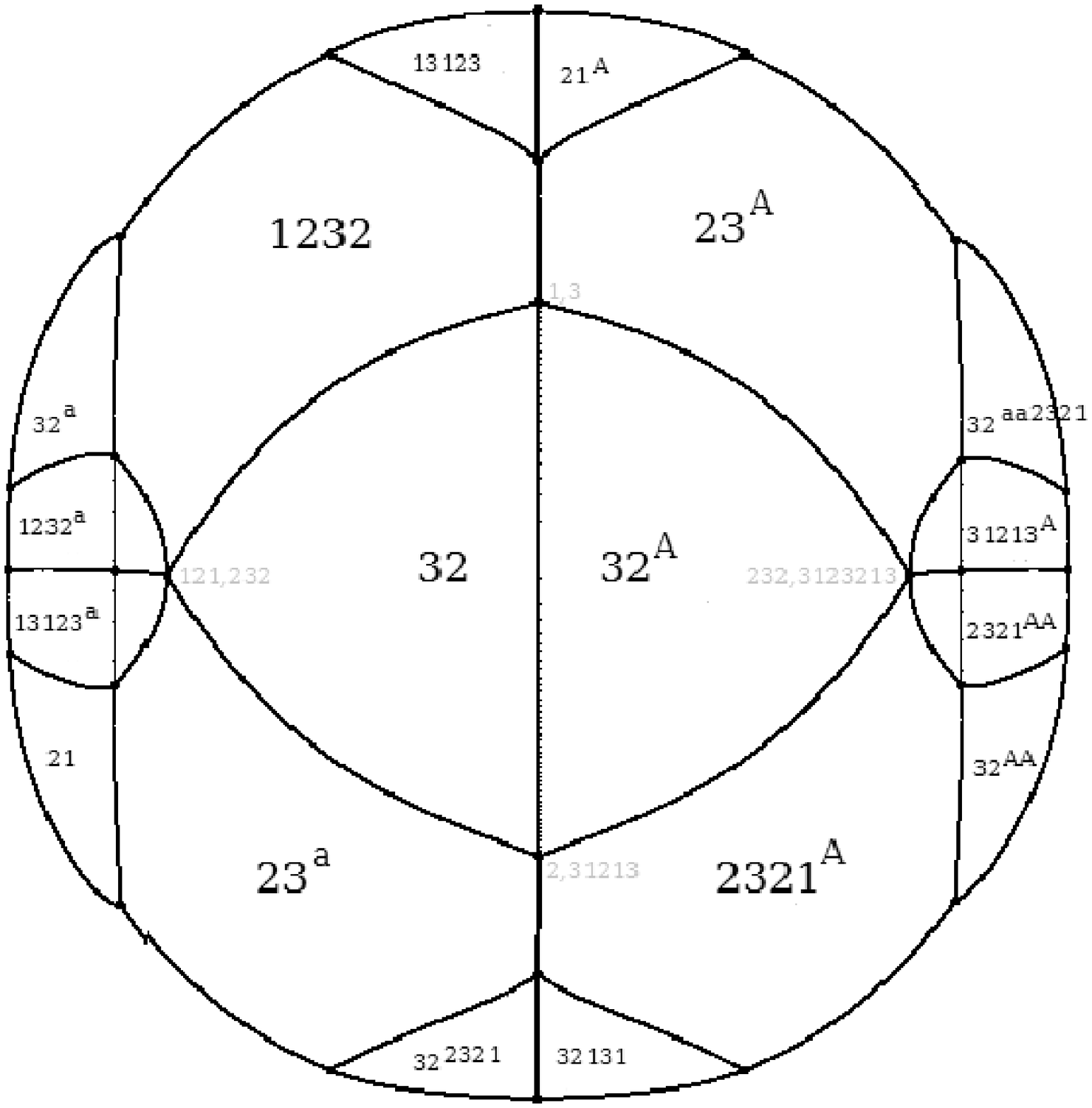, width=0.45\textwidth}}\hfill\,\\
\hfill\subfigure[2321]{\epsfig{figure=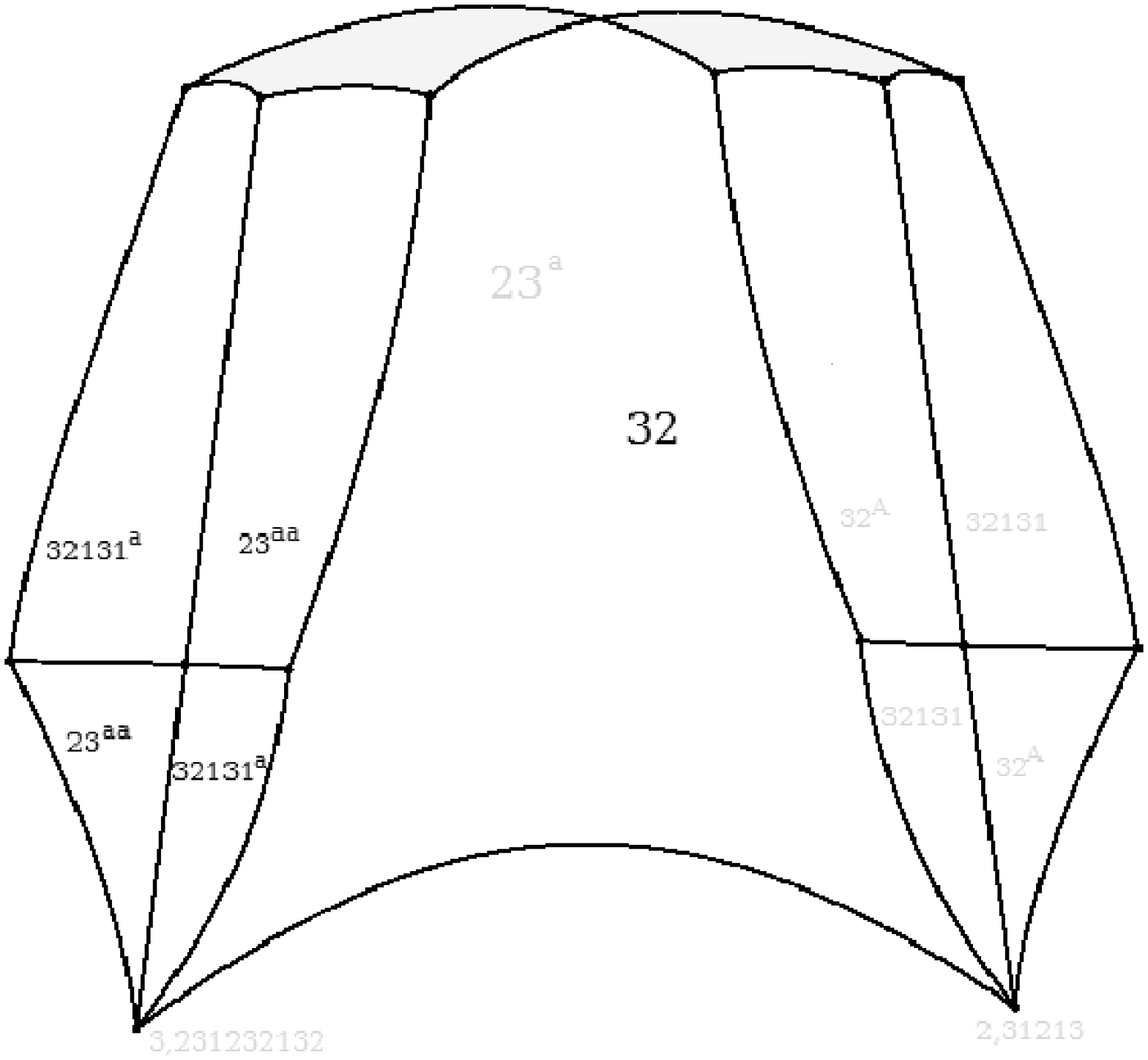, width=0.45\textwidth}}\hfill
\subfigure[1232]{\epsfig{figure=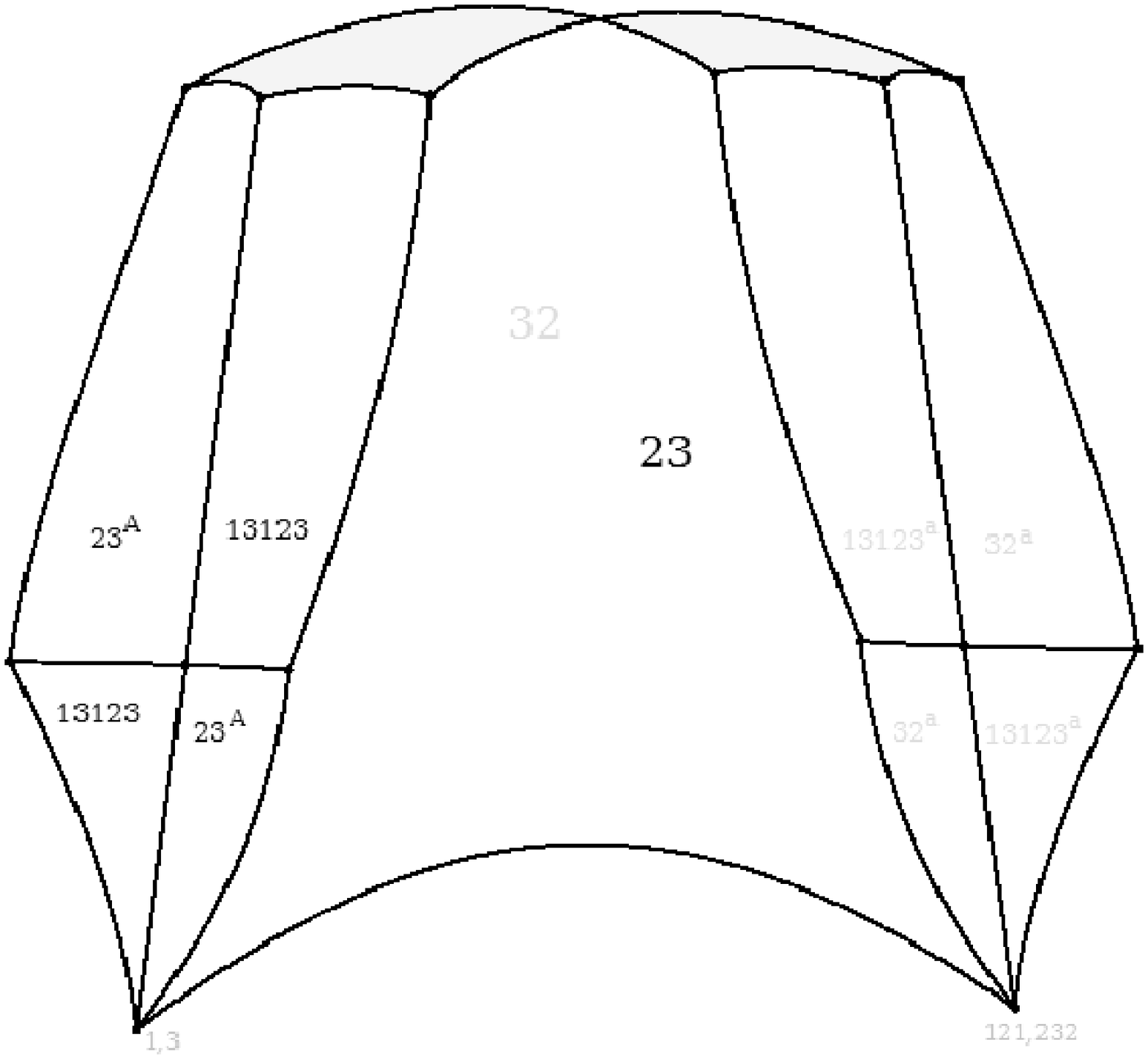, width=0.45\textwidth}}\hfill\,\\
\caption{Combinatorics of some faces of the Ford domain. In the faces
  for 23 and 32, there is an extra 16-gon that lies on the boundary at
  infinity, not drawn in the picture. For other faces, the boundary
  2-face is shaded in gray. The Ford domain has infinitely many faces,
  but every face is the image of one of the depicted 10 faces under an
  appropriate power of $a=2313$ (we write $A=3132$ for its inverse);
  see also Figure~\ref{fig:faces-2}.}\label{fig:faces-1}
\end{figure}

\begin{figure}[htbp]
\hfill\subfigure[12]{\epsfig{figure=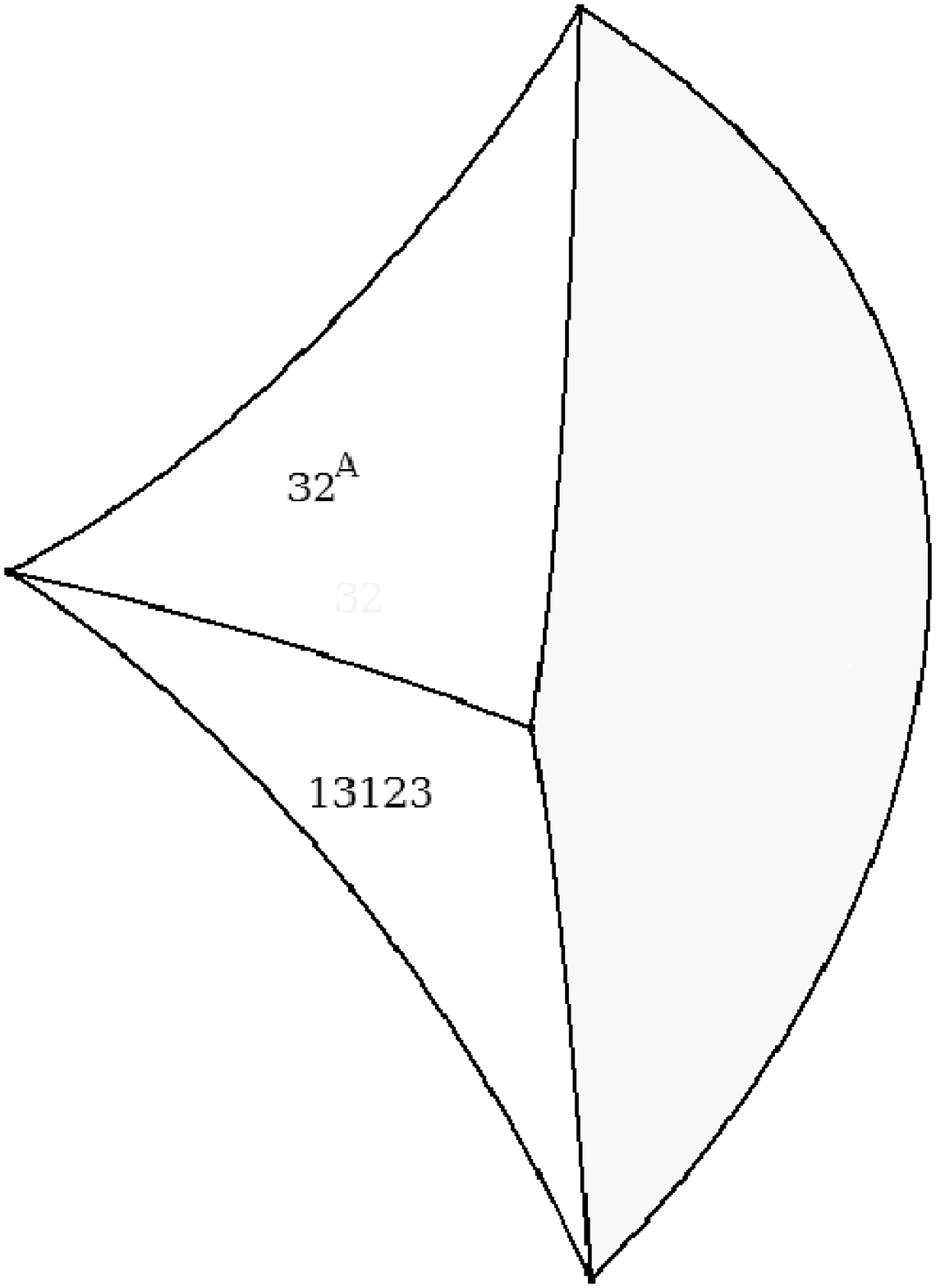, width=0.2\textwidth}}\hfill
\subfigure[21]{\epsfig{figure=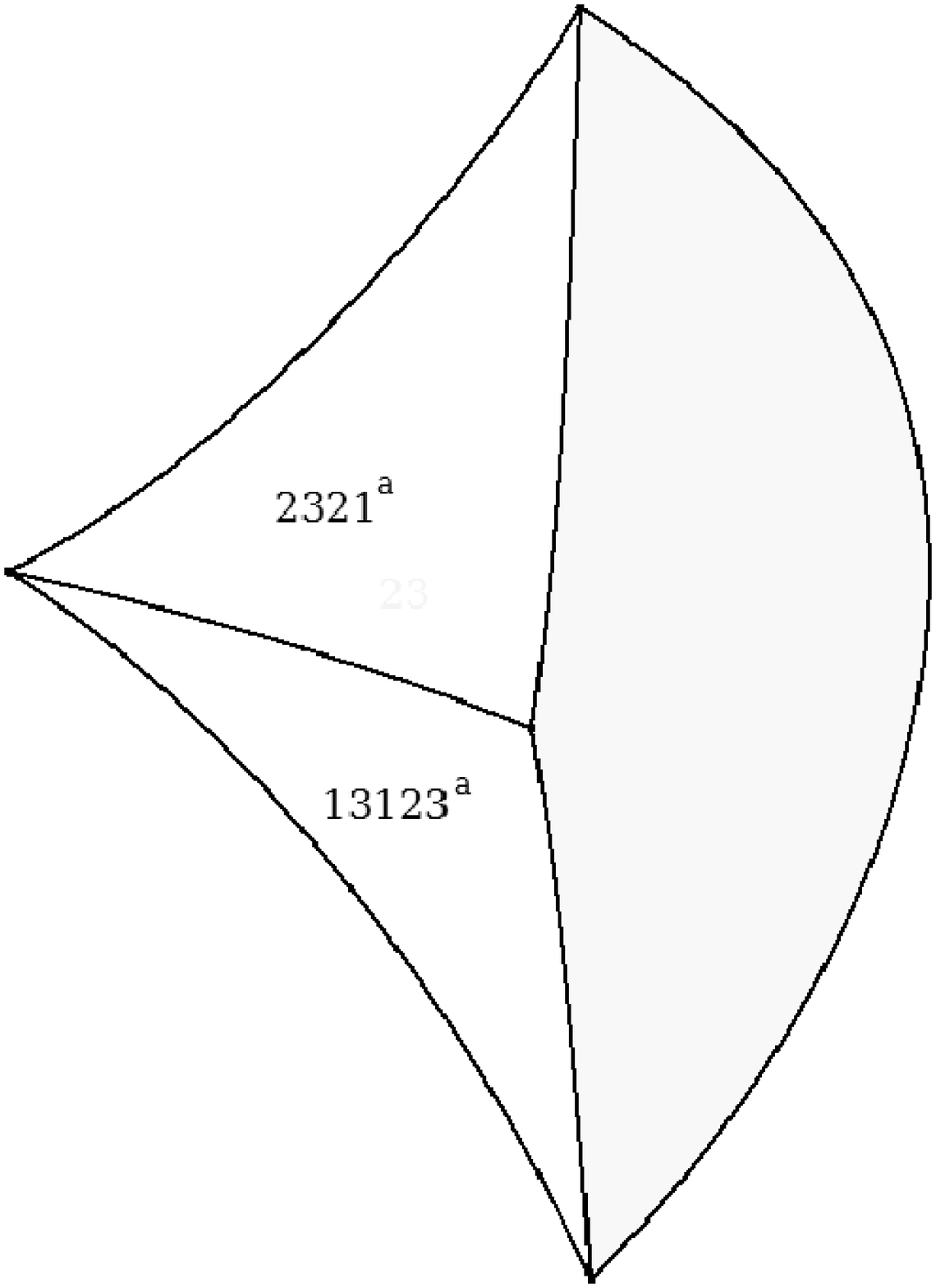, width=0.2\textwidth}}\hfill
\hfill\subfigure[32131232]{\epsfig{figure=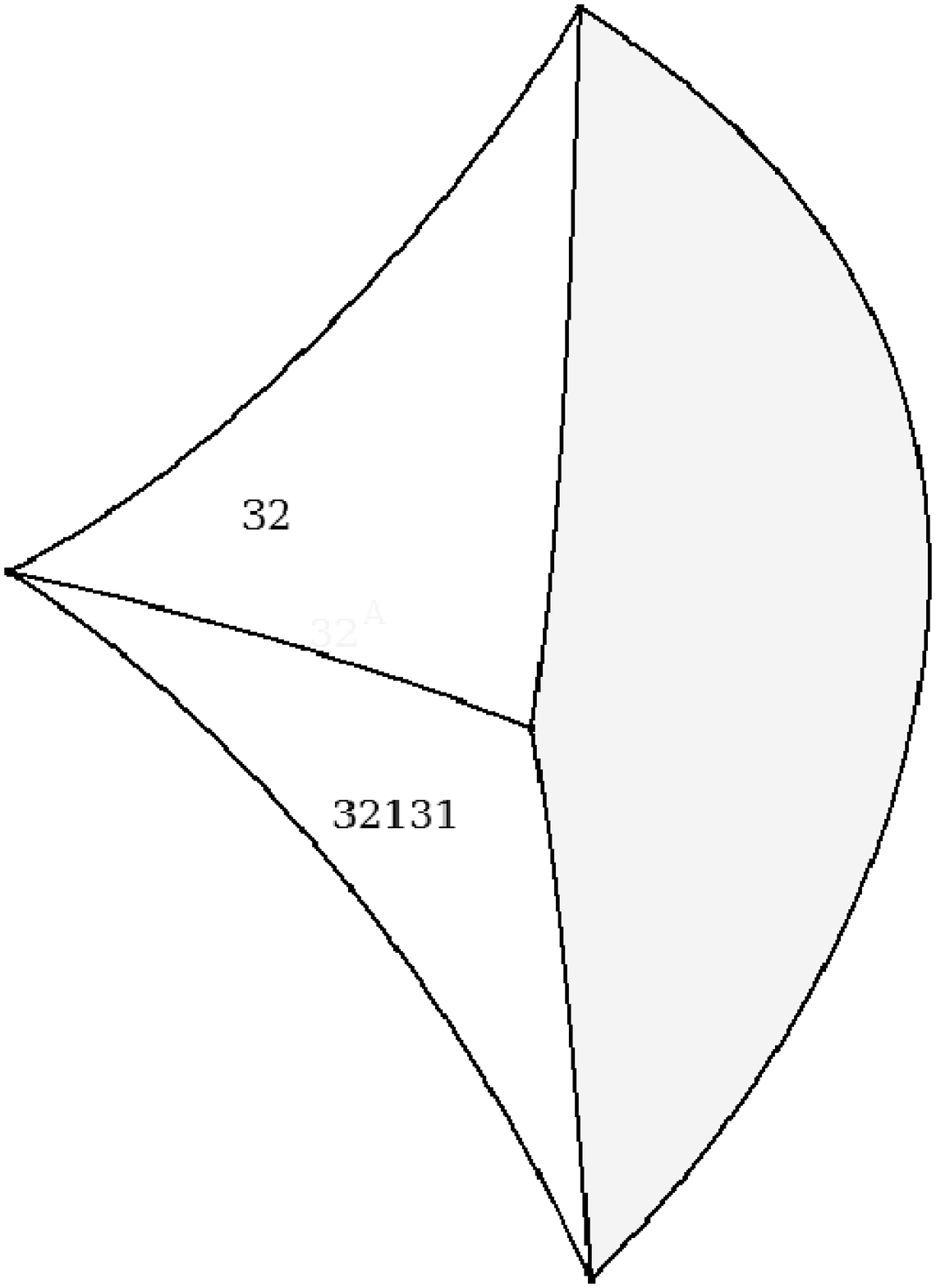, width=0.2\textwidth}}\hfill
\subfigure[23213123]{\epsfig{figure=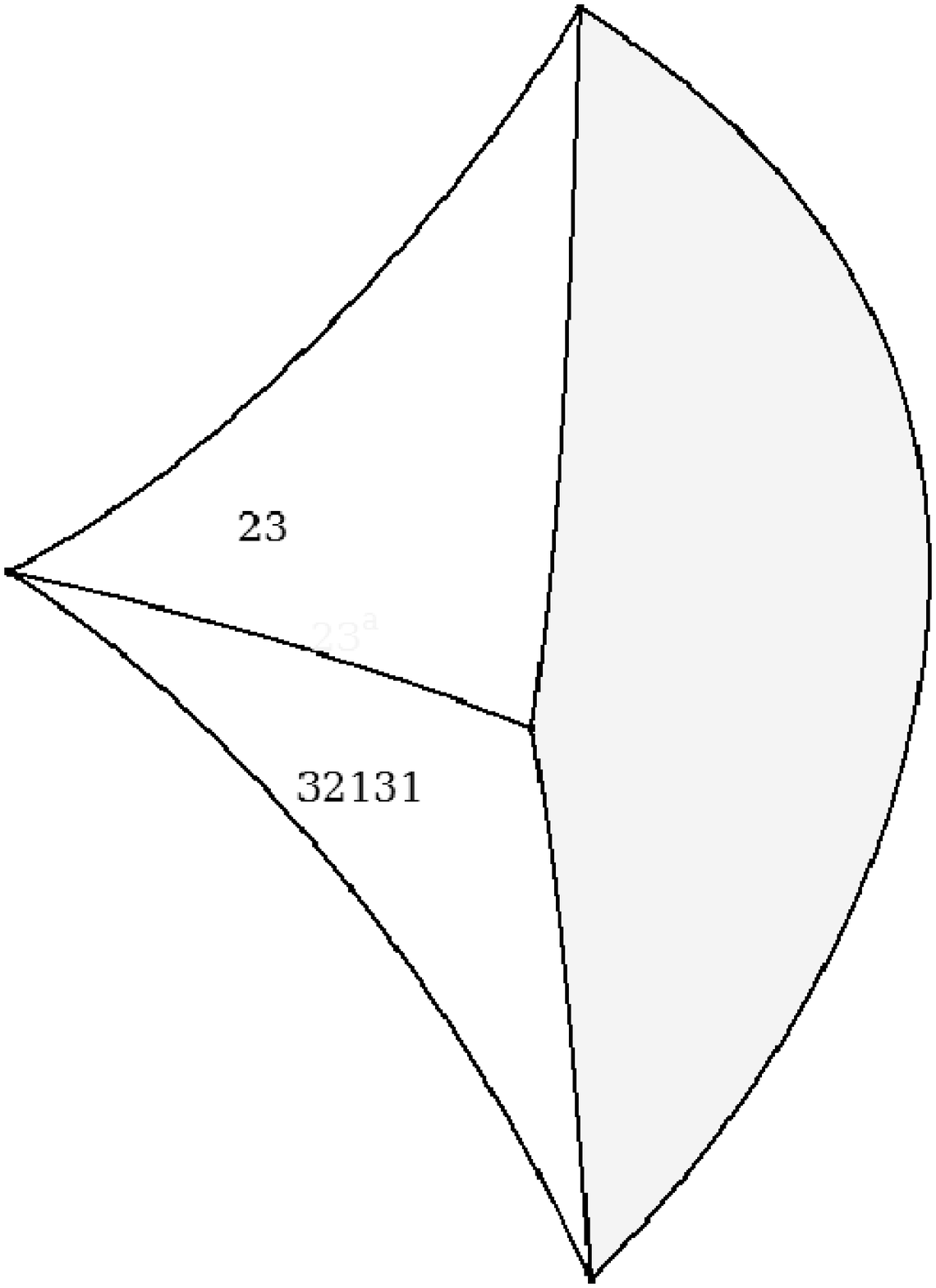, width=0.2\textwidth}}\hfill\,\\
\hfill\subfigure[232131]{\epsfig{figure=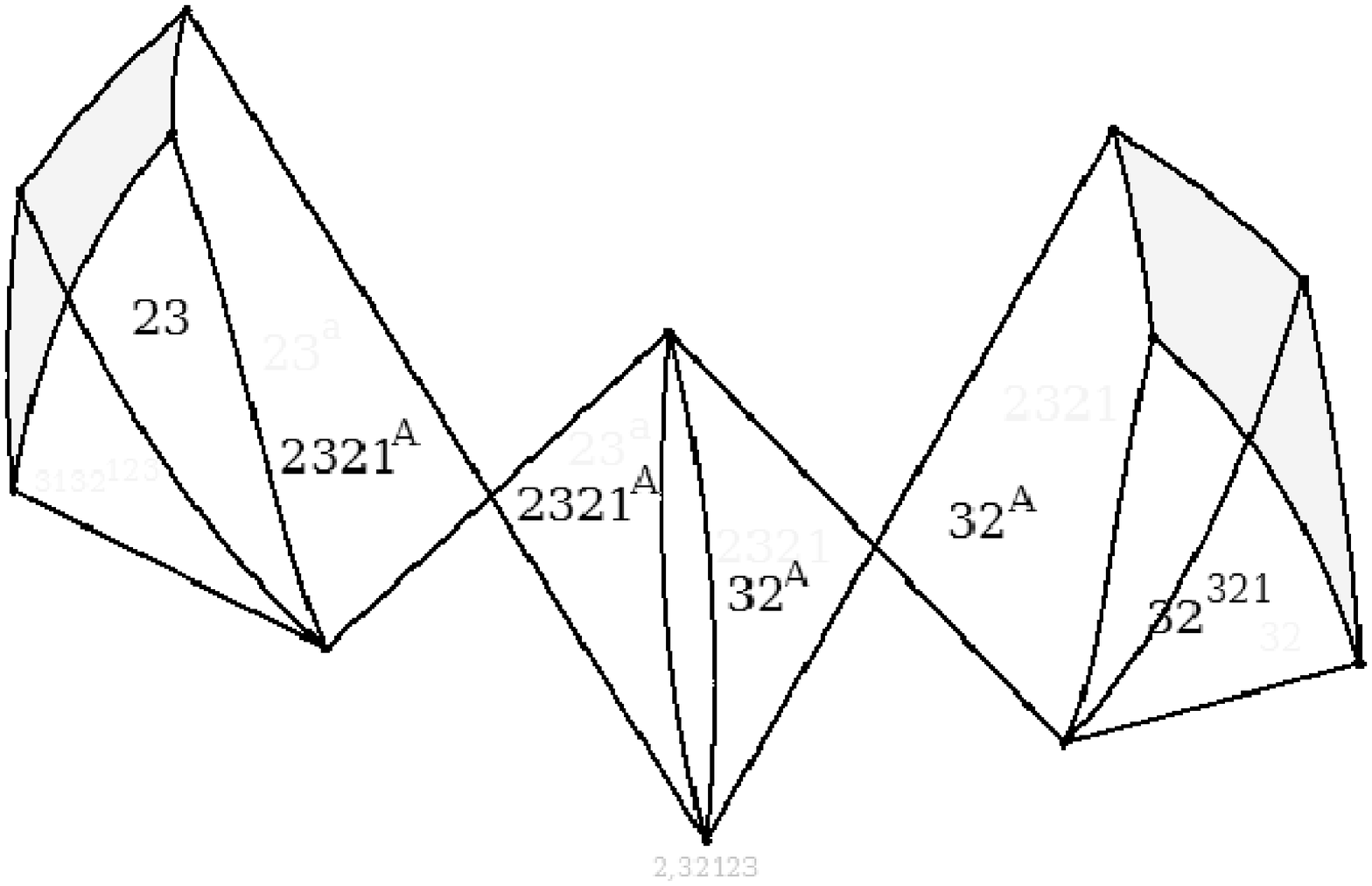, width=0.45\textwidth}}\hfill
\subfigure[131232]{\epsfig{figure=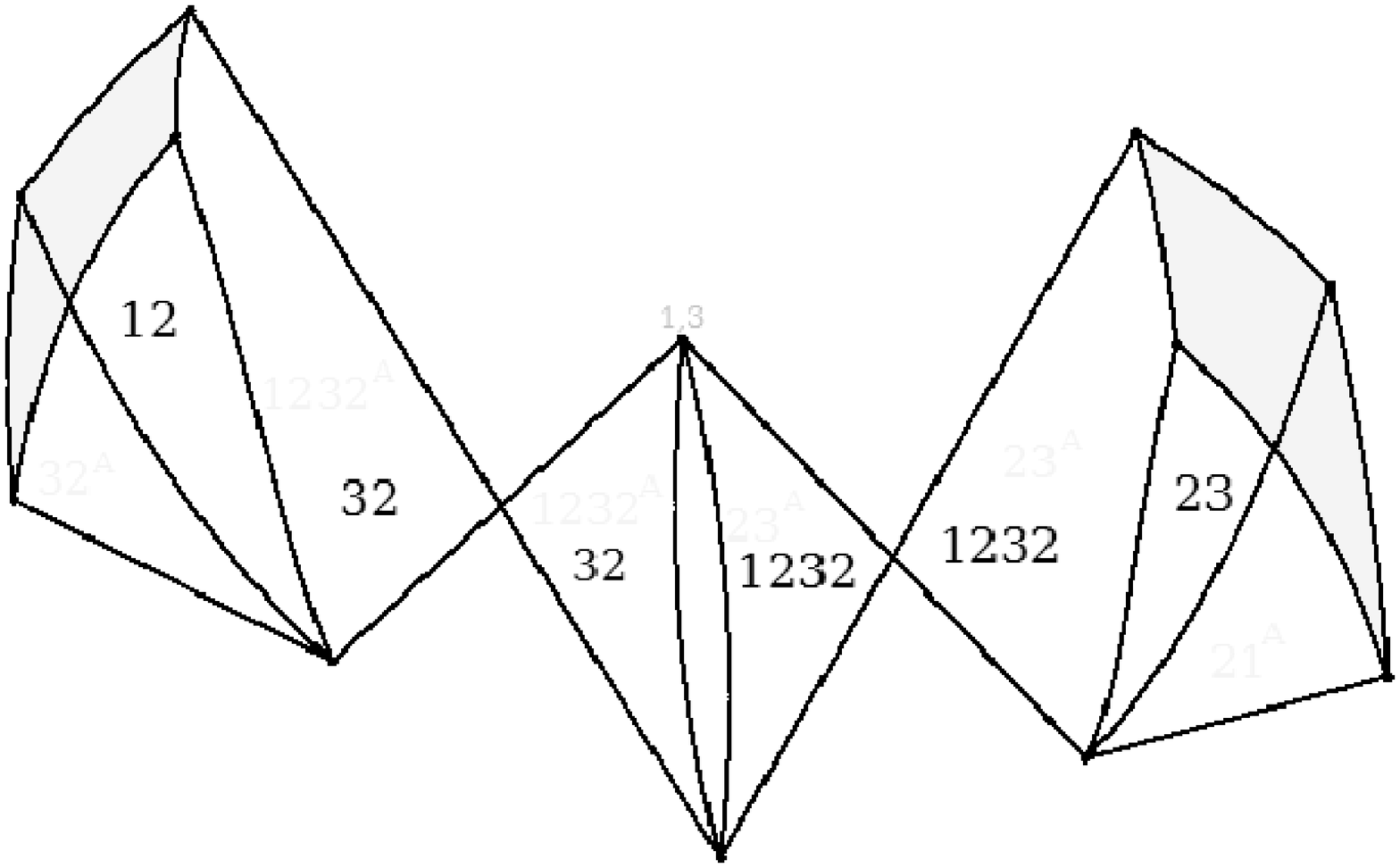, width=0.45\textwidth}}\hfill\,
\caption{More faces of the Ford domain.}\label{fig:faces-2}
\end{figure}

The corresponding ten polytopes $b_{g_1},\dots,b_{g_{10}}$ are depicted in
Figures~\ref{fig:faces-1} and~\ref{fig:faces-2}.
The pictures were obtained by parametrizing the 1-skeleton, and
mapping the 1-skeleton by a transformation that identifies the ambient
bisector with the unit ball in $\R^3$, in such a way that the
$\C$-slices are given by horizontal disks, and $\R$-slices are given
by vertical lines containing the $z$-axis. For an explanation of how
this can be done explicitly, see~\cite{derauxfalbel} for instance.

A priori, it is not clear why each face should have only finitely many
neighboring faces; this follows easily from the fact that the action
of $a$ in the $x$-coordinate is a translation, and the fact that the
spinal spheres in a Ford domain are given in Heisenberg coordinates by
bounded sets (note that if two spinal spheres are disjoint, then the
corresponding bisectors are disjoint).\label{page:disjoint}

In what follows, we write $F$ for $F_{\Gamma,p_\infty}$, and $E$ for
$\partial_\infty F\setminus \{p_\infty\}\subset \C\times\R$. We use Heisenberg
coordinates $(x,y,t)$, where $z=x+iy$. One easily computes the action
of $a=I_2I_3I_1I_3$ to be given by
\begin{equation}\label{eq:tsl}
(z,t)\mapsto(z+1,t-Im(z)).
\end{equation}
In particular, this map preserves every horizontal line in the plane
$Im(z)=0$. Each of these lines is a $\R$-circle (going through the
point $p_\infty$), i.e. it bounds a totally geodesic copies of
$H^2_\R$ in $H^2_\C$. The union of these real planes is the
\emph{invariant fan} of $a$, see~\cite{goldmanparker}.

For the sake of brevity, we write $p$ rather than $p_\infty$.
\begin{thm}\label{thm:poincare}
 \begin{enumerate}
   \item The maps $g_k$ define side pairings of $F$. More precisely,
     if $k$ is odd, $g_k(b_k)=b_{k+1}$ and $F\cap g_k(F)=b_{k+1}$;
     when $k$ is even, $g_k(b_k)=b_{k-1}$ and $F\cap g_k(F)=b_{k-1}$.
   \item These pairings satisfies the hypotheses of the Poincar{\'e}
     polyhedron theorem for cosets of the unipotent cyclic group
     $\langle a\rangle$.
   \item The group $\Gamma$ has the following presentation,
     $$
     \langle x_1,x_2 | x_1^3,x_2^3, (x_1x_2)^5 \rangle,
     $$
     where we have written $x_1=I_1I_2$, $x_2=I_2I_3$.
 \end{enumerate}
\end{thm}

\begin{pf}
(1) The key for checking this is to certify that the combinatorics given
  in Figures~\ref{fig:faces-1} and~\ref{fig:faces-2} are correct. We
  do not expand on the details, but this can be done because the
  entries of the generators, as well as the center of the Ford domain,
  can be given by entries a number field of small degree (here the
  degree is 4, see page~\pageref{pagenbfield}).

  One easily verifies that the isometries given in
  Table~\ref{tab:corefaces} define side pairings of the Ford domain,
  by computing several orbits under appropriate group
  elements. Clearly it is enough to work on the core faces, i.e. the
  representatives given in Table~\ref{tab:corefaces}. We will give
  some detail only for the first two faces, the other ones being
  entirely similar.

  For instance, the fact that $I_2I_3$ maps $b_1$ to $b_2$ follows
  from the fact that $I_2I_3$ does what is announced in
  Table~\ref{tab:orbits}, where we use the numbering of
  equation~\ref{eq:numbering}.
  \begin{table}
    \begin{tabular}{cccccccccc}
      Point & \#2 & \#3  & \#4  & \#5  & \#7   & \#8   & \#9 & \#11 & \#12\\
      Image & \#1 & \#22 & \#12 & \#11 & \# 26 & \# 14 & \#8 & \#23 & \#21\\[0.3cm]
      Point & \#13 & \#15 & \#17 & \#18 & \#19 & \#21 & \#22 & \#24 & \#32\\
      Image & \#27 & \#7  & \#43 & \#10 & \#41 & \# 4 & \#6  & \#18 & \#30
    \end{tabular}
    \caption{Table of correspondence of the 2-faces of $b_1$ and
      $b_2$, under the natural side pairing map $I_2I_3:B_1\rightarrow
      B_2$.}\label{tab:orbits}
  \end{table}
  Clearly by construction $I_2I_3(B_1)=B_2$, because
  $g_2=g_1^{-1}$. The first column of the table means that
  $I_2I_3(B_1\cap B_2)=B_1\cap B_2$, and this follows readily from the
  fact that $I_1I_2$ has order 3.

  The next column in the table says that $I_2I_3(B_1\cap B_3)=B_2\cap
  B_{22}$; this follows from the fact that
  $$
  I_2I_3(g_3p)=g_{22}p=a^{-1}g_2p.
  $$ 
  Equivalently, we claim
  $$
  23\cdot 2321 p=3132\cdot 23 p.
  $$ 
  This is an obvious consequence of the fact that $(23)^3=id$.

  In fact all other claims in the table are all consequences of the
  relations $(12)^3,(23)^3$ and $(31)^5$, as well as the fact that
  2313 fixes $p$.
  In any event, it should be clear that the claims in the table can
  readily be checked with a computer (at worst, one performs
  computation in the relevant number field).

(2) The second item is checked by tracing the Poincar{\'e}
  cycles. Recall that a \emph{ridge} is by definition a codimension
  two facet of $F$. Note that no ridge of the Ford domain is totally
  geodesic (this requires a computation, it amounts to saying that for
  $k\neq l$, $p$, $g_kp$ and $g_lp$ are never in a common complex
  line, or in other words, any choice of homogeneous coordinates for
  these three vectors produces a basis of $\C^3$). By the discussion on
  page~\pageref{page:disjoint}, only finitely many checks need to be made,
  since for $m$ large enough, $B_k\cap a^m B_l=\emptyset$.

   The ridges of $F$ are so-called Giraud disks, which are generic
   intersections of two bisectors; because the complex spines all
   intersect in $p$, we can think of the intersections as being
   coequidistant, and in particular their intersections are all smooth
   disks, equidistant of three points
   $p,g_kp,g_lp$, with $k\neq l$.

   Because of Giraud's theorem
   (see~\cite{giraud},~\cite{goldman},~\cite{deraux4445}), the ridges
   of $F$ are on precisely three bisectors, so the local tiling
   condition near generic ridges is actually a consequence of the
   existence of side-pairings.

(3) The explicit cycles are obtained by computing orbits of these
   triples of points under the side pairings; whenever a ridge in the
   cycle differs from the starting ridge by a power of $a$, we close
   the cycle up by that power of $a$ (see~\cite{derauxparkerpaupert}
   or~\cite{parkerbook}).

   We work out a few cycles, the other ones being similar. The ridge
   $b_1\cap b_2$ is sent to itself by $I_2I_3$, in fact
   $$
   p \stackrel{23}{\longrightarrow} 23(p) \stackrel{23}{\longrightarrow} 2323(p)=32(p) \stackrel{23}{\longrightarrow} p.
   $$
   This clearly gives a cycle transformation of order 3 preserving
   that ridge, so we get the relation
   $$
   (23)^3=id.
   $$

   The ridge $b_1\cap b_3$ is slightly more interesting. One checks
   (most conveniently with a computer), that
   $$
   p,g_1p,g_3p
   \stackrel{2313}{\longleftarrow}
   p,g_{21}p,g_{23}p
   \stackrel{321313}{\longleftarrow}
   g_{24}p,g_{21}p,p
   \stackrel{313213}{\longleftarrow}
   g_{2}p,p,g_{22}p   
   \stackrel{23}{\longleftarrow}
   p,g_1p,g_3p.
   $$
   The corresponding relation is
   $$
   2313\cdot 321313\cdot 313213 \cdot 23 = id,
   $$
   which can be simplified, using $232=323$, to
$$
   (12)^3=id.
$$
   One easily checks that $b_1\cap b_7$ gives
   $$
   (13)^5=id.
   $$
   Using the relations $(12)^3=(23)^3=(31)^5=id$, one checks that the
   other cycle relations can be reduced to a trivial relation.    
\end{pf}

\begin{prop}
  \begin{enumerate}
    \item The only elliptic elements in $\Gamma$ are conjugates of powers of
  $12$, $23$ or $31$; in particular, no elliptic element of $\Gamma$
  fixes any point in $\partial_\infty H^2_\C$. 
    \item The only parabolic
  elements in $\Gamma$ are conjugates of powers of $2313$.
  \end{enumerate}
\end{prop}

\begin{pf}
(1) As mentioned in the proof of Theorem~\ref{thm:poincare}, the ridge
cycles are all conjugate to powers of $12$, $23$ or $31$. One then
considers cycles of lower-dimensional facets, namely 1-faces and
vertices. 

The cycles of 1-faces turn out to be trivial, and the only non-trivial
vertex cycles correspond to the fixed points of $13$ and that of
$321323$ (or conjugates of these under some power of $a$).

(2) Since ideal vertices all have trivial stabilizers, the only parabolic
elements in the group are the ones stabilizing the center of the Ford
domain, which is by construction a fundamental domain modulo cosets of
$\langle 2313\rangle$
\end{pf}

The combinatorial structure of $\partial E$ can be gathered from the
shaded 2-faces in Figures~\ref{fig:faces-1} and~\ref{fig:faces-2}
(apart from the first two faces, where the corresponding boundary
14-gon is not shown on the picture). It may seem somewhat reminiscent
of the boundary of the real hyperbolic Ford domain, but it is quite
different (unlike the case of the spherical CR uniformization of the
figure eight knot complement, see~\cite{derauxdeformfig8}).

Because of the shearing by the imaginary part of $z$ in
formula~\ref{eq:tsl}, it is not that easy to produce a meaningful
2-dimensional picture of $\partial E$, which is topologically a
cylinder. Note that the $x$-axis is entirely outside $E$, and it gives
a core curve for a solid cylinder (in $\partial_\infty H^2_\C$, one
gets a solid torus pinched at $p_\infty$). This means that $E$ is the
complement of a topological solid cylinder; it is in fact a horotube,
in Schwartz's terminology~\cite{richBook}.

The determination of the topology of the manifold at infinity is
somewhat delicate. From the combinatorial description of the $\partial
E$ (together with the action of the cyclic group generated by 2313),
one can compute the fundamental group of the manifold. Indeed, one can
start with a presentation of the fundamental group of the 1-skeleton,
and include a relation saying that each loop coming from the boundary
of a 2-cell becomes trivial in the 2-skeleton. 

The bookkeeping of this computation is of course prohibitingly lengthy
when performed by hand, but it is fairly well suited to calculations
performed by the computer. The end result is that the manifold at
infinity does have the same fundamental group as \verb|m009|, and one
can check the peripheral subgroups are preserved under this
isomorphism. 

%\bibliographystyle{plain}
%\bibliography{biblio}

\end{document}